\documentstyle[12pt,amssymb]{article}

\newcommand{\rest}{\restriction}

\makeatletter
\let\DOTSI\relax
\def\RIfM@{\relax\ifmmode}%
\def\FN@{\futurelet\next}%
\newcount\intno@
\def\iint{\DOTSI\intno@\tw@\FN@\ints@}%
\def\iiint{\DOTSI\intno@\thr@@\FN@\ints@}%
\def\iiiint{\DOTSI\intno@4 \FN@\ints@}%
\def\idotsint{\DOTSI\intno@\z@\FN@\ints@}%
\def\ints@{\findlimits@\ints@@}%
\newif\iflimtoken@
\newif\iflimits@
\def\findlimits@{\limtoken@true\ifx\next\limits\limits@true
 \else\ifx\next\nolimits\limits@false\else
 \limtoken@false\ifx\ilimits@\nolimits\limits@false\else
 \ifinner\limits@false\else\limits@true\fi\fi\fi\fi}%
\def\multint@{\int\ifnum\intno@=\z@\intdots@                                %1
 \else\intkern@\fi                                                          %2
 \ifnum\intno@>\tw@\int\intkern@\fi                                         %3
 \ifnum\intno@>\thr@@\int\intkern@\fi                                       %4
 \int}%                                                                     %5
\def\multintlimits@{\intop\ifnum\intno@=\z@\intdots@\else\intkern@\fi
 \ifnum\intno@>\tw@\intop\intkern@\fi
 \ifnum\intno@>\thr@@\intop\intkern@\fi\intop}%
\def\intic@{\mathchoice{\hskip.5em}{\hskip.4em}{\hskip.4em}{\hskip.4em}}%
\def\negintic@{\mathchoice
 {\hskip-.5em}{\hskip-.4em}{\hskip-.4em}{\hskip-.4em}}%
\def\ints@@{\iflimtoken@                                                    %1
 \def\ints@@@{\iflimits@\negintic@\mathop{\intic@\multintlimits@}\limits    %2
  \else\multint@\nolimits\fi                                                %3
  \eat@}%                                                                   %4
 \else                                                                      %5
 \def\ints@@@{\iflimits@\negintic@
  \mathop{\intic@\multintlimits@}\limits\else
  \multint@\nolimits\fi}\fi\ints@@@}%
\def\intkern@{\mathchoice{\!\!\!}{\!\!}{\!\!}{\!\!}}%
\def\plaincdots@{\mathinner{\cdotp\cdotp\cdotp}}%
\def\intdots@{\mathchoice{\plaincdots@}%
 {{\cdotp}\mkern1.5mu{\cdotp}\mkern1.5mu{\cdotp}}%
 {{\cdotp}\mkern1mu{\cdotp}\mkern1mu{\cdotp}}%
 {{\cdotp}\mkern1mu{\cdotp}\mkern1mu{\cdotp}}}%
\def\Let@{\relax\iffalse{\fi\let\\=\cr\iffalse}\fi}%
\def\vspace@{\def\vspace##1{\crcr\noalign{\vskip##1\relax}}}%
\def\multilimits@{\bgroup\vspace@\Let@
 \baselineskip\fontdimen10 \scriptfont\tw@
 \advance\baselineskip\fontdimen12 \scriptfont\tw@
 \lineskip\thr@@\fontdimen8 \scriptfont\thr@@
 \lineskiplimit\lineskip
 \vbox\bgroup\ialign\bgroup\hfil$\m@th\scriptstyle{##}$\hfil\crcr}%
\def\Sb{_\multilimits@}%
\def\endSb{\crcr\egroup\egroup\egroup}%
\def\Sp{^\multilimits@}%

\newdimen\ex@
\def\rightarrowfill@#1{$#1\m@th\mathord-\mkern-6mu\cleaders
 \hbox{$#1\mkern-2mu\mathord-\mkern-2mu$}\hfill
 \mkern-6mu\mathord\rightarrow$}%
\def\leftarrowfill@#1{$#1\m@th\mathord\leftarrow\mkern-6mu\cleaders
 \hbox{$#1\mkern-2mu\mathord-\mkern-2mu$}\hfill\mkern-6mu\mathord-$}%
\def\leftrightarrowfill@#1{$#1\m@th\mathord\leftarrow\mkern-6mu\cleaders
 \hbox{$#1\mkern-2mu\mathord-\mkern-2mu$}\hfill
 \mkern-6mu\mathord\rightarrow$}%
\def\overrightarrow{\mathpalette\overrightarrow@}%
\def\overrightarrow@#1#2{\vbox{\ialign{##\crcr\rightarrowfill@#1\crcr
 \noalign{\kern-\ex@\nointerlineskip}$\m@th\hfil#1#2\hfil$\crcr}}}%

\def\overleftarrow{\mathpalette\overleftarrow@}%
\def\overleftarrow@#1#2{\vbox{\ialign{##\crcr\leftarrowfill@#1\crcr
 \noalign{\kern-\ex@\nointerlineskip}$\m@th\hfil#1#2\hfil$\crcr}}}%
\def\overleftrightarrow{\mathpalette\overleftrightarrow@}%
\def\overleftrightarrow@#1#2{\vbox{\ialign{##\crcr\leftrightarrowfill@#1\crcr
 \noalign{\kern-\ex@\nointerlineskip}$\m@th\hfil#1#2\hfil$\crcr}}}%
\def\underrightarrow{\mathpalette\underrightarrow@}%
\def\underrightarrow@#1#2{\vtop{\ialign{##\crcr$\m@th\hfil#1#2\hfil$\crcr
 \noalign{\nointerlineskip}\rightarrowfill@#1\crcr}}}%

\def\underleftarrow{\mathpalette\underleftarrow@}%
\def\underleftarrow@#1#2{\vtop{\ialign{##\crcr$\m@th\hfil#1#2\hfil$\crcr
 \noalign{\nointerlineskip}\leftarrowfill@#1\crcr}}}%
\def\underleftrightarrow{\mathpalette\underleftrightarrow@}%
\def\underleftrightarrow@#1#2{\vtop{\ialign{##\crcr$\m@th\hfil#1#2\hfil$\crcr
 \noalign{\nointerlineskip}\leftrightarrowfill@#1\crcr}}}%
\newcount\GRAPHICSTYPE
\GRAPHICSTYPE=\z@
\def\GRAPHICSPS#1{%
 \ifnum\GRAPHICSTYPE=\@ne language "PS", include "#1"\else ps: #1\fi
 }%
\def\graffile#1#2#3#4{%
 \leavevmode\raise -#4 \hbox{%
  \raise #3 \hbox{\rule{0.003in}{0.003in}\special{#1}}%
  }%
 {\raise -#4 \hbox to #2 {\vrule height#3 width\z@ depth\z@\hfil}}%
 }%
\def\draftbox#1#2#3#4{%
 \leavevmode\raise -#4 \hbox{%
  \frame{\rlap{\protect\tiny #1}\hbox to #2%
   {\vrule height#3 width\z@ depth\z@\hfil}%
  }%
 }%
}%
\newcount\draft
\draft=\z@
\def\GRAPHIC#1#2#3#4#5{%
 \ifnum\draft=\@ne \draftbox{#2}{#3}{#4}{#5}%
  \else \graffile{#1}{#3}{#4}{#5}%
  \fi
 }%
\def\addtoLaTeXparams#1{\edef\LaTeXparams{\LaTeXparams #1}}%
\def\doFRAMEparams#1{\readFRAMEparams#1\end}%
\def\readFRAMEparams#1{%
 \ifx#1\end%
  \let\next=\relax
  \else
  \ifx#1i\dispkind=\z@\fi
  \ifx#1d\dispkind=\@ne\fi
  \ifx#1f\dispkind=\tw@\fi
  \ifx#1t\addtoLaTeXparams{t}\fi
  \ifx#1b\addtoLaTeXparams{b}\fi
  \ifx#1p\addtoLaTeXparams{p}\fi
  \ifx#1h\addtoLaTeXparams{h}\fi
  \let\next=\readFRAMEparams
  \fi
 \next
 }%
\def\IFRAME#1#2#3#4#5{\GRAPHIC{#5}{#4}{#1}{#2}{#3}}%
\def\DFRAME#1#2#3#4{%
 \begin{center}\GRAPHIC{#4}{#3}{#1}{#2}{\z@}\end{center}%
 }%
\def\FFRAME#1#2#3#4#5#6#7{%
 \begin{figure}[#1]%
  \begin{center}\GRAPHIC{#7}{#6}{#2}{#3}{\z@}\end{center}%
  \caption{\label{#5}#4}%
  \end{figure}%
 }%
\newcount\dispkind%
\def\FRAME#1#2#3#4#5#6#7#8{%
 \def\LaTeXparams{}%
 \dispkind=\z@
 \def\LaTeXparams{}%
 \doFRAMEparams{#1}%
 \ifnum\dispkind=\z@\IFRAME{#2}{#3}{#4}{#7}{#8}\else
  \ifnum\dispkind=\@ne\DFRAME{#2}{#3}{#7}{#8}\else
   \ifnum\dispkind=\tw@
    \edef\@tempa{\noexpand\FFRAME{\LaTeXparams}}%
    \@tempa{#2}{#3}{#5}{#6}{#7}{#8}%
    \fi
   \fi
  \fi
 }%
\def\limfunc#1{\hbox{\rm #1}}%
\long\def\QQQ#1#2{\long\expandafter\def\csname#1\endcsname{#2}}%
\def\QTP#1{}%
\long\def\QQA#1#2{}%
\def\QTR#1#2{{\csname#1\endcsname #2}}%(gp) Is this the best?
\long\def\TeXButton#1#2{#2}%
\def\EXPAND#1[#2]#3{}%
\def\NOEXPAND#1[#2]#3{}%
\def\LaTeXparent#1{}%
\def\QTagDef#1#2#3{}%
\def\QQfnmark#1{\footnotemark}

\@ifundefined{INDEX}{\def\INDEX#1#2{}{}}{}%
\@ifundefined{SUBINDEX}{\def\SUBINDEX#1#2#3{}{}{}}{}%
\def\initial#1{\bigbreak{\raggedright\large\bf #1}\kern 2\p@\penalty3000}%
\@ifundefined{abstract}{%
 \def\abstract{%
  \if@twocolumn
   \section*{Abstract (Not appropriate in this style!)}%
   \else \small 
   \begin{center}{\bf Abstract\vspace{-.5em}\vspace{\z@}}\end{center}%
   \quotation 
   \fi
  }%
 }{%
 }%
\@ifundefined{endabstract}{\def\endabstract
  {\if@twocolumn\else\endquotation\fi}}{}%
\@ifundefined{maketitle}{\def\maketitle#1{}}{}%
\@ifundefined{affiliation}{\def\affiliation#1{}}{}%
\@ifundefined{proof}{\def\proof{\paragraph{Proof. }}}{}%
\@ifundefined{endproof}{}{}%
\@ifundefined{newfield}{\def\newfield#1#2{}}{}%
\@ifundefined{chapter}{\def\chapter#1{\par(Chapter head:)#1\par }%
 \newcount\c@chapter}{}%
\@ifundefined{part}{\def\part#1{\par(Part head:)#1\par }}{}%
\@ifundefined{section}{\def\section#1{\par(Section head:)#1\par }}{}%
\@ifundefined{subsection}{\def\subsection#1%
 {\par(Subsection head:)#1\par }}{}%
\@ifundefined{subsubsection}{\def\subsubsection#1%
 {\par(Subsubsection head:)#1\par }}{}%
\@ifundefined{paragraph}{\def\paragraph#1%
 {\par(Subsubsubsection head:)#1\par }}{}%
\@ifundefined{subparagraph}{\def\subparagraph#1%
 {\par(Subsubsubsubsection head:)#1\par }}{}%
\@ifundefined{therefore}{}{}%
\@ifundefined{backepsilon}{}{}%
\@ifundefined{yen}{}{}%
\@ifundefined{registered}{\def\registered{\relax\ifmmode{}\r@gistered
                                                \else$\m@th\r@gistered$\fi}%
 \def\r@gistered{^{\ooalign
  {\hfil\raise.07ex\hbox{$\scriptstyle\rm\hbox{R}$}\hfil\crcr
  \mathhexbox20D}}}}{}%
\@ifundefined{Eth}{}{}%
\@ifundefined{eth}{}{}%
\@ifundefined{Thorn}{}{}%
\@ifundefined{thorn}{}{}%
\@ifundefined{degree}{}{}%
\def\BibTeX{{\rm B\kern-.05em{\sc i\kern-.025em b}\kern-.08em
    T\kern-.1667em\lower.7ex\hbox{E}\kern-.125emX}}%
\newdimen\theight
\def\Column{%
 \vadjust{\setbox\z@=\hbox{\scriptsize\quad\quad tcol}%
  \theight=\ht\z@\advance\theight by \dp\z@\advance\theight by \lineskip
  \kern -\theight \vbox to \theight{%
   \rightline{\rlap{\box\z@}}%
   \vss
   }%
  }%
 }%
\def\qed{%
 \ifhmode\unskip\nobreak\fi\ifmmode\ifinner\else\hskip5\p@\fi\fi
 \hbox{\hskip5\p@\vrule width4\p@ height6\p@ depth1.5\p@\hskip\p@}%
 }%
\def\miss{\hbox{\vrule height2\p@ width 2\p@ depth\z@}}%
\def\tcol#1{{\baselineskip=6\p@ \vcenter{#1}} \Column}  %
\makeatother

\newtheorem{theorem}{Theorem}
\newtheorem{lemma}[theorem]{Lemma}
\newtheorem{proposition}[theorem]{Proposition}

\begin{document}

\author{Paul C. Eklof\\
University of California, Irvine  \and Matthew Foreman\\
%EndAName
University of California, Irvine  \and Saharon Shelah\\
Hebrew University and Rutgers University }
\thanks{%
Thanks to Rutgers University for its support of this research through its
funding of the first and third authors' visits to Rutgers.} 
\thanks{%
The second author thanks the NSF, Grant No. DMS-9203726, for partial support.}
\thanks{%
Partially supported by Basic Research Fund, Israeli Academy of Sciences.
Pub. No. 520 }   %EndAName
\title{On invariants for $\omega _1$-separable groups}
\maketitle

\begin{abstract}
We study the classification of $\omega_1$-separable groups using
Ehrenfeucht-Fra\"\i ss\'e games and prove a strong classification result
assuming PFA, and a strong non-structure theorem assuming $\diamondsuit$.
\end{abstract}

\section*{Introduction}

An $\omega _1$-separable (or $\aleph _1$-separable) group is an abelian
group such that every countable subset is contained in a free direct summand
of the group. In particular, therefore, an $\omega _1$-separable group is $%
\aleph _1$-free, i.e., every countable subgroup is free. The structure of $%
\omega _1$-separable groups of cardinality $\aleph _1$ was investigated in 
\cite{E} and \cite{M}; most of the results proved there required
set-theoretic assumptions beyond ZFC. (See also \cite[Chap. VIII]{EM} for an
exposition of these results.) Specifically, assuming Martin's Axiom (MA)
plus $\neg $CH or the stronger Proper Forcing Axiom (PFA), one can prove
nice structure and classification results; these results are not theorems of
ZFC since counterexamples exist assuming CH or ``prediction principles''
like $\diamondsuit $. In \cite[Remark 3.3]{E} it is asserted that a
construction given there under the assumption of CH (or even $2_{}^{\aleph
_0}<2^{\aleph _1}$) of two non-isomorphic $\omega _1$-separable groups

\begin{quote}
``is strong evidence for the claim that in a model of CH there is no
possible meaningful classification of all $\omega _1$-separable groups. It
is difficult to see what conceivable scheme of classification could
distinguish between [the groups constructed here].''
\end{quote}

\noindent  But, in fact, the Helsinki school of model theory provides a
scheme for distinguishing between such groups. It is our aim here to use the
methodology of the Helsinki school --- which involves Ehrenfeucht-Fra\"\i
ss\'e games (cf. \cite{MSV}, \cite{O} or \cite{V}) --- to strengthen the dichotomy referred to above: that is, to
obtain strong classification results assuming PFA, and a strong
``non-structure theorem'' assuming $\diamondsuit .$

We begin by describing the Ehrenfeucht-Fra\"\i ss\'e (or EF) games, after
which we can state our results more precisely. If $\alpha $ is an ordinal
and $A$ and $B$ are any structures, the game $EF_\alpha (A,B)$ is played
between two players $\forall $ and $\exists $ who take turns choosing
elements of $A\cup B$ through $\alpha $ rounds. Specifically, in each round $%
\forall $ picks first an element of either $A$ or $B$; and then $\exists $
picks an element of the other structure. The result is, at the end, two
sequences $(a_\nu )_{\nu <\alpha }$ and $(b_\nu )_{\nu <\alpha }$ of
elements of, respectively, $A$ and $B$. Player $\exists $ wins if and only
if the function $f$ which takes $a_\nu $ to $b_\nu $ is a partial
isomorphism; otherwise $\forall $ wins. If $A$ and $B$ have cardinality $%
\kappa $, $\exists $ has a winning\ strategy for $EF_\kappa (A,B)$ if and
only if $A$ and $B$ are isomorphic. (Let $\forall $ list all the elements of 
$A\cup B$ during his moves.)

We consider variations of these games defined using trees. Given any tree $T$%
, we define the game $EF(A,B;T)$: the game is played as before except that
player $\forall $ must also, whenever it is his turn, pick a node of the
tree strictly above his previous choices (thus his successive choices will
form a branch --- a linearly ordered subset --- of the tree). The game ends
when $\forall $ can no longer pick a node above his previous choices; the
criterion for winning is as before, that is, $\exists $ wins if and only if
the function $f$ defined by the play is a partial isomorphism. We write $%
A\equiv ^TB$ if $\exists $ has a winning strategy in the game $EF(A,B;T)$.
For the purposes of motivation consider first the case $\alpha =\omega $.
(Our interest is in the case $\alpha =\omega _1$.) In this case, we consider
only well-founded trees, i.e., trees without infinite branches; then for
every such $T$, each play of the game $EF(A,B;T)$ is finite. (So 
$EF(A,B;T)$ may be regarded as an approximation to the game 
$EF_\omega (A,B)$.) Scott's Theorem
implies that for each countable $A$ there is a countable ordinal $\beta $
such that if $T_\beta $ is any tree of rank $\beta $, then for any countable 
$B$, $B$ is isomorphic to $A$ if and only if $A\equiv ^{T_\beta ^{}}B$ . In
terms of infinitary languages, $A$ is determined up to isomorphism (among
countable structures) by a sentence of $L_{\infty \omega }$ of rank $\beta $.

For structures of cardinality $\aleph _1$, it is natural to look at
approximations to the EF game of length $\omega _1$ and use trees which may
have countably infinite branches, but do not have branches of cardinality $%
\aleph _1$; we call these {\it bounded trees}. For such $T$, each play of
the game $EF(A,B;T)$ will end after countably many moves. We will say $A$\
is $T${\it -equivalent} to $B$ if $A\equiv ^TB$. This relation provides 
a possible way of distinguishing between the $\omega_1$-separable groups 
constructed in \cite{E} under the assumption of CH (cf. the remark 
after the quotation 
above).

 By a theorem of Hyttinen 
\cite{Hy}, the entire class of bounded trees determines $A$ up to
isomorphism; that is, if $A$ and $B$ are of cardinality $\aleph _1$ and $%
A\equiv ^TB$ for all bounded trees, then $A$ is isomorphic to $B$. The 
structure of the class of bounded trees is much more complicated than 
that of the class of well-founded trees (cf. \cite{V}). However,
in contrast to the situation for countable structures, there is not always a
single tree which suffices to describe $A$ up to isomorphism. Specifically,
Hyttinen and Tuuri \cite{HyT} proved (assuming CH) that there is a linear
order $A$ of cardinality $\aleph _1$ such that for every bounded tree $T$
there is a linear order $B_T$ of cardinality $\aleph _1$ such that $A\equiv
^TB_T$ but $A$ is not isomorphic to $B_T$. They call this result a {\it %
non-structure theorem }for $A$. It can be translated in terms of infinitary
languages and says that there is no complete description of $A$ in a certain
strong language $M_{\omega _2\omega _1}$ (which we shall not define here).

A similar non-structure theorem for $p$-groups was proved by Mekler and
Oikkonen \cite{MO}; their theorem is proved by carrying over to $p$-groups,
by means of a Hahn power construction, the result of Hyttinen and Tuuri.
Whether the analogous result for $\aleph _1$-free groups is a theorem of ZFC
+ CH  remains open, but when we consider the question for $\aleph _1$%
-separable groups, we obtain an independence result, which is the subject of
this paper. In the first section we prove (with the help of the structural
results referred to above) that assuming PFA

\begin{quote}
{\it if }$A${\it \ and }$B${\it \ are }$\omega _1${\it -separable groups of
cardinality }$\aleph _1${\it \ such that }$A\equiv ^{\omega
^2+\omega }B${\it , then they are isomorphic (where }$\omega ^2+\omega ${\it %
\ is the countable ordinal regarded as a --- linearly ordered --- tree).}
\end{quote}

\noindent See Theorem \ref{stform}. Thus a single, simple, tree contains
enough information to classify any $\omega _1$-separable group --- in the
precise sense that a single sentence of $M_{\omega _2\omega _1}$ of ``tree
rank'' $\omega ^2+\omega $ completely describes $A$.

In section 2 we show, assuming $\diamondsuit $, that not only 
does $\omega ^2 + \omega$ not have the property above, but
 for {\em any} bounded tree $T$,
 there are non-isomorphic $\omega _1$-separable groups $A^T$ and $B^T$
 of cardinality $\aleph_1$ which
cannot be separated by $T$, in the sense that $A^T\equiv ^T B^T$. (See Theorem 
\ref{onetree}.) The construction in section 2 is strengthened in section 3
to obtain a non-structure theorem (Theorem \ref{nonstr}.):

\begin{quote}
{\it there is an $\omega _1$-separable group }$A${\it \ of cardinality }$%
\aleph _1${\it \ such that for every bounded tree }$T${\it \ there is an $%
\omega _1$-separable group }$B^T${\it \ of cardinality }$\aleph _1${\it \
which is not isomorphic to }$A${\it \ but is }$T$-{\it equivalent} to $A$.
\end{quote}

\noindent (Note that $A$ does not depend on $T$.)

We shall make use, at times, of the following simple lemma, where $A^{*}$
denotes the dual of $A$, i.e. $\limfunc{Hom}(A,{\Bbb Z)}$.

\begin{lemma}
\label{stein} Suppose $A\subseteq B$ and $A^{\prime }\subseteq
B^{\prime }\subseteq C^{\prime }$ where $C^{\prime }/B^{\prime }$ is $\aleph
_1$-free, $B/A$ is countable and $(B/A)^{*}=0$. If $\theta :B\rightarrow
C^{\prime }$ such that $\theta [A]\subseteq A^{\prime }$, then $\theta
[B]\subseteq B^{\prime }$.
\end{lemma}

\TeXButton{Proof}{\proof} $\theta $ induces a homomorphism: $B/A\rightarrow
C^{\prime }/A^{\prime }$. By the hypotheses, the composition of this map
with the canonical surjection: $C^{\prime }/A^{\prime }\rightarrow C^{\prime
}/B^{\prime }$ must be zero; that is, $\theta [B]\subseteq B^{\prime }$. $%
\Box $

\section{A structure theorem}

An $\aleph _1$-separable group $A$ of cardinality $\aleph _1$ is
characterized by the property that it has a {\it 
filtration},
 that is, a
continuous chain $\{A_\nu :\nu <\omega _1\}$ of countable free subgroups
whose union is $A$ and is such that $A_0=0$ and for all $\nu $, $A_{\nu +1}$
is a direct summand of $A$. We say that two $\aleph _1$-separable groups $A$
and $B$ are {\it quotient-equivalent} if and only if they have filtrations, $%
\{A_\nu :\nu <\omega _1\}$ and $\{B_\nu :\nu <\omega _1\}$, respectively,
such that for every $\alpha <\omega _1$,  $A_{\alpha +1}/A_\alpha $ is  
isomorphic to $A_{\alpha +1}^{\prime }/A_\alpha
^{\prime }$. We say that  $A$ and $B$ are {\it %
filtration-equivalent} if and only if they satisfy the stronger condition
that  for every $\alpha <\omega _1$ there is a {\it level-preserving
isomorphism} $\theta _\alpha :A_{\alpha +1}\rightarrow B_{\alpha +1}$, i.e.,
an isomorphism such that for every $\nu \leq \alpha $, $\theta [A_\nu
]=B_\nu $. Under the assumption of MA + $\neg $CH, filtration-equivalence
implies isomorphism.

In \cite{M} (see also \cite[Chap. VIII]{EM}) it is proved under the
hypothesis of the Proper Forcing Axiom, PFA, that $\aleph _1$-separable groups of
cardinality $\aleph _1$ have a nice structure theory. More precisely, it 
is shown that, under PFA, every $\aleph _1$-separable group of
cardinality $\aleph _1$ is {\em in standard form}. (Roughly, this means 
that they have a ``classical" construction. We will give a definition 
below.)  Our goal in this
section is to use that theory to prove the following:

\begin{theorem}
\label{pfa}\ (PFA) $\omega ^2+\omega $ is a universal equivalence tree for
the class of $\aleph _1$-separable abelian groups of cardinality $\aleph _1$%
. That is, any two $\aleph _1$-separable abelian groups of cardinality $%
\aleph _1$ which are $\omega ^2+\omega $-equivalent are isomorphic.
\end{theorem}

We shall see in the next section that this is not a theorem of ZFC. We begin
with a weaker result.

\begin{proposition}
\label{quoteq} If $A$ and $A^{\prime }$ are strongly $\aleph _1$-free groups
of cardinality $\aleph _1$ which are $\omega 2$-equivalent, then they are
quotient equivalent.
\end{proposition}

\TeXButton{Proof}{\proof} Suppose that $\tau $ is a w.s. for $\exists $. Let 
$C$ be a cub such that if $\alpha \in C$ then for any $n\in \omega $, as
long as the first $n$ moves of $\forall $ are in $A_\alpha \cup A_\alpha
^{\prime }$, the replying moves of $\exists $ given by $\tau $ are also in $%
A_\alpha \cup A_\alpha ^{\prime }$. If $A$ and $A^{\prime }$ are not
quotient-equivalent, there exists $\alpha \in C$ such that $A_{\alpha
+1}/A_\alpha \oplus {\Bbb Z}^{(\omega )}$ is not isomorphic to $A_{\alpha
+1}^{\prime }/A_\alpha ^{\prime }\oplus {\Bbb Z}^{(\omega )}$. Now let $%
\forall $ play the game so that during the first $\omega $ moves he makes
sure that all elements of $A_\alpha \cup A_\alpha ^{\prime }$ are played;
the result, since $\tau $ is a w.s., is that an isomorphism $f:A_\alpha
\rightarrow A_\alpha ^{\prime }$ is obtained.

Then in the next $\omega $ moves, $\forall $ plays so that all, and only,
the elements of $A_\beta \cup A_\beta ^{\prime }$ are played for some $\beta
\geq \alpha +1$. This is possible by using a bijection of $\omega $ with $%
\omega \times \omega $. The result is an extension of $f$ to an isomorphism $%
f^{\prime }:A_\beta \rightarrow A_\beta ^{\prime }.$ Then, since $A_\beta
/A_{\alpha +1}$ and A$_\beta ^{\prime }/A_{\alpha +1}^{\prime }$ are free,
we have $A_\beta /A_\alpha \cong A_{\alpha +1}/A_\alpha \oplus A_\beta
/A_{\alpha +1}$ and similarly on the other side. Since $f^{\prime }$ induces
an isomorphism of $A_\beta /A_\alpha $ with $A_\beta ^{\prime }/A_\alpha
^{\prime }$, we obtain a contradiction of the choice of $\alpha $. $\Box $

\bigskip\ 

Suppose $A$ is an $\aleph_1$-separable group of cardinality $\aleph_1$ 
with a filtration  $\{A_\nu \colon \nu
\in \omega _1\}$, and let $E=\{\delta :A_\delta $ is not a direct
summand of $A\}$; $A$ is said to be {\em in standard form} if:

(1) it has a coherent system of projections $\{\pi _\nu 
\colon \nu \notin E\}$, i.e., projections $\pi _\nu { }\colon A \rightarrow  A_\nu { }$
 with 
the property that for all  $\nu  < \tau $  in $\omega _1 \setminus  E$,  
$\pi _\nu { } \circ  \pi _\tau  = \pi _\nu { }$; and

(2) for every $\delta \in
E$ there is a ladder $\eta _\delta $ on $\delta $ and a subset $Y_\delta $
of $A_{\delta +1}$ such that $A_{\delta +1}=A_\delta +\langle Y_\delta
\rangle $ and

\begin{quote}
($\dagger $) for all $y\in \langle Y_\delta \rangle $ and all $\nu <\delta $
with $\nu \notin E$, 
 $\pi _\nu (y)=\sum_{\alpha \in S}(\pi _{\alpha
+1}(y)-\pi _\alpha (y))$ where $S = \{ \alpha \in \limfunc{rge}(\eta 
_\delta) \colon \alpha < \nu \}$. 
\end{quote}

\noindent (Here a ladder on $\delta $ means a strictly increasing function $\eta
_\delta :\omega \rightarrow \delta $ with $\limfunc{rge}(\eta _\delta
)\subseteq \omega _1\setminus E$ and $\sup \limfunc{rge}(\eta _\delta
)=\delta $.) This property is actually stronger than the usual definition of
standard form (because of the assertion about the ladder); it can be shown
that the Proper Forcing Axiom (PFA) implies that every strongly $\aleph _1$%
-free group of cardinality $\aleph _1$ has this property (by essentially the
same proof as in \cite[Thm. VIII.3.3]{EM}).

Let  $K_\alpha = \ker (\pi _\alpha)$ and let  
$K_{\alpha, \alpha + 1}  = K_\alpha \cap  A_{\alpha + 1} $. Notice that we can replace any $y$ in $Y_\delta $ by $y+u$ where $u\in
K_{\alpha ,\alpha +1}$ for some $\alpha \in \delta \setminus E$, and we will
still have a generating set of $A_{\delta +1}$ over $A_\delta $ which
satisfies ($\dagger $). Also we can, and will, assume that $A_{\nu +1}/A_\nu 
$ has infinite rank for every $\nu \notin E$.

\begin{lemma}
\label{gens} Suppose $A$ is in standard form. Then there is a filtration $%
\{A_\nu :\nu \in \omega _1\}$ of $A$ and for each $\delta \in E=\{\delta
:A_\delta $ is not a direct summand of $A\}$, there are: a ladder $\eta
_\delta $ on $\delta $; and a subset $\bar y_\delta =\{y_{\delta ,i}:i\in I\}
$ of $A_{\delta +1}$ which is linearly independent mod $A_\delta $ such that
if $\beta _n=\eta _\delta (n)$:

\begin{enumerate}
\item  for all $n\in \omega $, $\beta _n\notin E$; and

\item  $A_{\delta +1}$ is generated mod $A_\delta $ by a set of elements of
the form 
\begin{equation}
\label{geners}\frac{t(\bar y_\delta )-a}d
\end{equation}
where $t(\bar y_\delta )$ is a linear combination of the elements of $\bar
y_\delta $, $d\in {\Bbb Z}$, and $a\in \oplus _{n\in \omega }K_{\beta
_n,\beta _n+1}$.
\end{enumerate}

\noindent Moreover, given $\mu <\delta $, we can choose $\eta _\delta $ such
that $\eta _\delta (0)>\mu $.
\end{lemma}

\TeXButton{Proof}{\proof} Let $Y_\delta $ and $\eta _\delta $ be as in the
definition of standard form above. Let $\bar y_\delta =\{y_{\delta ,i}:i\in
I\}$ be a maximal linearly independent subset of $Y_\delta $. By the remark
preceding the lemma we can (by replacing $y_{\delta ,i}$ by $y_{\delta ,i}+u$
for some $u$) assume that $\eta _\delta (0)>\mu $.

If $d$ divides $t(\bar y_\delta )$ mod $A_{\nu +1}$ for some integer $d$ and
linear combination $t(\bar y_\delta )$, then $d$ divides $t(\bar y_\delta
)-a $ where $a=\pi _{\nu +1}($ $t(\bar y_\delta ))=\sum_{\beta \in S}\pi
_{\beta ,\beta +1}($ $t(\bar y_\delta ))$ for some finite subset $S\subseteq 
\limfunc{rge}(\eta _\delta )$. $\Box $

\bigskip\ 

\begin{proposition}
\label{fileq} Let $G$ and $G^{\prime }$ be $\aleph _1$-separable groups such
that $G$ is in standard form. Suppose that they have filtrations $\{G_\nu
:\nu \in \omega _1\}$ and $\{G_\nu ^{\prime }:\nu \in \omega _1\}$
respectively such that the filtration of $G$ attests that $G$ is in standard
form and $E=\{\nu \in \omega _1:G_\nu $ is not a summand of $G\}=\{\nu \in
\omega _1:G_\nu ^{\prime }$ is not a summand of $G^{\prime }\}$. Suppose
also that for all limit ordinals $\delta $, given a ladder $\eta _\delta $
on $\delta $, there is an isomorphism $\theta _\delta :G_{\delta
+1}\rightarrow G_{\delta +1}^{\prime }$ such that for all $n\in \omega $, $%
\theta _\delta [G_{\eta _\delta (n)}]=G_{\eta _\delta (n)}^{\prime }$ and $%
\theta _\delta [G_{\eta _\delta (n)+1}]=G_{\eta _\delta (n)+1}^{\prime }$.
Then $G$ and $G^{\prime }$ are filtration-equivalent.
\end{proposition}

\TeXButton{Proof}{\proof} We can assume that the filtration of $G$ is as in
Lemma \ref{gens}. We prove by induction on $\nu $ the following:

\begin{quote}
if $\mu <\nu $ and $\mu ,\nu \in \omega _1\setminus E$ and $f:G_\mu
\rightarrow G_\mu ^{\prime }$ is a level-preserving isomorphism, then $f$
extends to a level-preserving isomorphism $g:$ $G_\nu \rightarrow G_\nu
^{\prime }$.
\end{quote}

If $\nu =\tau +1$ where $\tau \notin E$, then the result follows easily by
induction and the fact that $G_\nu /G_\tau $ and $G_\nu ^{\prime }/G_\tau
^{\prime }$ are free. If $\nu $ is a limit ordinal, choose a ladder $\zeta
_\nu $ on $\nu $ such that $\zeta _\nu (0)>\mu $ and for all $n$, $\zeta _\nu
(n)\notin E$, and extend $f$ successively, by induction, to $g_n:G_{\zeta
_\nu (n)}\rightarrow G_{\zeta _\nu (n)}^{\prime }$, and let $g=\cup _ng_n$.

The crucial case is when $\nu =\delta +1$ where $\delta \in E$. Let $\eta
_\delta $ be as in Lemma \ref{gens} with $\eta _\delta (0)>\mu ,$ and let $%
\theta _\delta $ be the corresponding isomorphism given by the hypothesis of
this Proposition. Let $C_{\delta ,n}=K_{\beta _n,\beta _n+1}$. By induction,
extend $f$ to a level-preserving isomorphism $f_0:G_{\eta _\delta
(0)}\rightarrow G_{\eta _\delta (0)}^{\prime }$ and then extend it to $%
g_0:G_{\eta _\delta (0)+1}\rightarrow G_{\eta _\delta (0)+1}^{\prime }$ by
letting $g_0 \rest C_{\delta ,0}=\theta _\delta  \rest C_{\delta ,0}$. Clearly $g_0$ is
level-preserving. By induction extend $g_0$ to a level-preserving $%
f_1:G_{\eta _\delta (1)}\rightarrow G_{\eta _\delta (1)}^{\prime }$ and then
to $g_1:G_{\eta _\delta (1)+1}\rightarrow G_{\eta _\delta (1)+1}^{\prime }$
by letting $g_1 \rest C_{\delta ,1}=\theta _\delta  \rest C_{\delta ,1}$. Continuing in
this way we obtain level-preserving isomorphisms $g_n:G_{\eta _\delta
(n)+1}\rightarrow G_{\eta _\delta (n)+1}^{\prime }$ for each $n$. Let $%
\tilde g=\cup _ng_n:G_\delta \rightarrow G_\delta ^{\prime }$.

By Lemma 4, $G_{\delta +1}$ is generated mod $G_\delta $ by a set of elements of the
form 
$$
\frac{t(\bar y_\delta )-a}d 
$$
where $a\in \oplus _{n\in \omega }C_{\delta ,n}$; hence $G_{\delta
+1}^{\prime }$ is generated mod $G_\delta ^{\prime }$ by elements 
$$
\frac{t(\theta _\delta (\bar y_\delta ))-\theta _\delta (a)}d\hbox{.} 
$$
But then since $\tilde g(a)=\theta _\delta (a)$ for each such $a$ by
construction, we can extend $\tilde g$ to $g:G_{\delta +1}\rightarrow
G_{\delta +1}^{\prime }$ by sending each $y_{\delta ,i}$ in $\bar y_\delta $
to $\theta _\delta (y_{\delta ,i})$. Since $\bar y_\delta $ is linearly
independent over $G_\delta $ this is a well-defined homomorphism. $\Box $

\bigskip\ 

\begin{theorem}
{\sc \label{stform}\ }Suppose $A$ and $A^{\prime }$ are $\aleph _1$%
-separable groups of cardinality $\aleph _1$ and at least one of them is in
standard form. If $A$ and $A^{\prime }$ are $\omega ^2+\omega $-equivalent,
then they are filtration-equivalent.
\end{theorem}

\TeXButton{Proof}{\proof} We can suppose that $A$ is in standard form, and
that we have chosen a filtration, $\{A_\nu \colon \nu \in \omega _1\}$ which
attests to that fact. Moreover, we can assume that if $\delta \in E=\{\delta
:A_\delta $ is not a direct summand of $A\},$ then $(A_{\delta +1}/A_\delta
)^{*}=0$. (Use Stein's Lemma \cite[Exer. 3, p. 112]{EM}, and replace $%
A_{\delta +1}$ by a direct summand, if necessary.)

Since $A$ is quotient-equivalent to $A^{\prime }$ by Proposition \ref{quoteq}
, we can assume that there is a filtration $\{A_\nu ^{\prime }\colon \nu \in
\omega _1\}$ of $A^{\prime }$ such that $E=\{\delta :A_\delta ^{\prime }$ is
not a direct summand of $A^{\prime }\}$ and for $\delta \in E$, $A_{\delta
+1}^{\prime }/A_\delta ^{\prime }\cong $ $A_{\delta +1}/A_\delta ^{}=0$, so
in particular $(A_{\delta +1}^{\prime }/A_\delta ^{\prime })^{*}=0 $.

Fix a bijection $\psi _{\alpha \beta }:\omega \rightarrow (A_\beta \setminus
A_\alpha )\cup (A_\beta ^{\prime }\setminus A_\alpha ^{\prime })$ for each $%
\alpha <\beta $. Let $\psi =\{\psi _{\alpha \beta }:\alpha <\beta <\omega
_1\}$.

Whenever we talk about moves in a game, we refer to the game EF$_{\omega
^2+\omega }(A,A^{\prime })$. Given a strictly increasing finite sequence of
countable ordinals $\alpha _1<\alpha _2<...<\alpha _n$, we will say that $%
\forall ${\it \ plays according to $\psi $ and }$\left\langle \alpha
_1,\alpha _2,...,\alpha _n\right\rangle $ {\it for the first }$\omega n${\it %
\ moves\ }if the $\omega k+\ell $ move of player $\forall $ is $\psi
_{\alpha _k\alpha _{k+1}}(\ell )$ for $k=0,...,n-1$ and $\ell \in \omega $
(where $\alpha _0=0$).

Suppose that $\tau $ is a w.s. for $\exists $ in the game EF$_{\omega
^2+\omega }(A,A^{\prime })$. Let $C$ be the set of all $\delta <\omega _1$
such that for any integers $n>0$ and $m\geq 0$ and any ordinals $\alpha
_1<\alpha _2<...<\alpha _n<\delta $, if $\forall $ plays according to $\psi $
and $\left\langle \alpha _1,\alpha _2,...,\alpha _n\right\rangle $ for the
first $\omega n$ moves and then plays any elements of $A_\delta $ for the
next $m$ moves, then the responses of $\exists $ using $\tau $ are all in $%
A_\delta \cup A_\delta ^{\prime }$.

Then $C$ is a cub: for the proof of unboundedness, note that there are only
countably many possibilities that one has to close under: choice of $n$ and $%
m$, choice of $\alpha _1<\alpha _2<...<\alpha _n$, and choice of moves $%
\omega n,\omega n+1,...\omega n+m-1$. (The earlier moves are determined by
the $\psi _{\alpha _k\alpha _{k+1}}$ and by $\tau $.)

There is a continuous strictly increasing function $\tilde h:\omega
_1\rightarrow \omega _1$ whose range is $C$. Define $h:\omega _1\rightarrow
\omega _1$ by%
$$
h(\beta )=\left\{ 
\begin{array}{ll}
\tilde h(\beta )+1 & \hbox{if }\beta \hbox{ is a successor and }\tilde
h(\beta )\in E \\ \tilde h(\beta ) & \hbox{otherwise} 
\end{array}
\right. 
$$
Let $G_\alpha =A_{h(\alpha )}$ and $G_\alpha ^{\prime }=A_{h(\alpha
)}^{\prime }$. Then for successor $\beta $, $G_\beta $ is a summand of $A$
and for limit $\delta $ $G_\delta =A_{\tilde h(\delta )}=\cup _{\beta
<\delta }A_{h(\beta )}=\cup _{\beta <\delta }G_\beta $, so $\{G_\alpha
:\alpha \in \omega _1\}$ (resp. $\{G_\alpha ^{\prime }:\alpha \in \omega
_1\} $) is a filtration of $A$ (resp. $A^{\prime }$). Given a limit ordinal $%
\delta $ and a ladder $\eta _\delta $ on $\delta $, it follows --- from Lemma 
\ref{stein} and the definition of $C$ --- that there is an isomorphism $%
\theta _\delta :G_{\delta +1}\rightarrow G_{\delta +1}^{\prime }$ such that
for all $n\in \omega $, $\theta _\delta [G_{\eta _\delta (n)}]=$ $G_{\eta
_\delta (n)}^{\prime }$ and $\theta _\delta [G_{\eta _\delta (n)+1}]=$ $%
G_{\eta _\delta (n)+1}^{\prime }$. In fact, $\theta _\delta $ is the partial
isomorphism which results because $\exists $ wins the game where the $\omega
k+\ell $ move of $\forall $ is 
$$
\psi _{ h(\eta _\delta (n)), h(\eta _\delta (n)+1)}(\ell ) 
$$
when $k = 2n $, and is
$$
\psi _{ h(\eta _\delta (n) + 1), h(\eta _\delta (n+1))}(\ell ) 
$$
when $k=2n+1$, and the $\omega ^2+m$ move of $\forall $ is $\psi
_{h(\delta ), h(\delta +1)}(m)$.

Thus we have satisfied the hypotheses of Proposition \ref{fileq} so we
conclude that $A$ and $A^{\prime }$ are filtration-equivalent. $\Box $

\bigskip 

Now we can prove Theorem \ref{pfa}. PFA implies that every strongly $\aleph
_1$-free abelian groups of cardinality $\aleph _1$ is $\aleph _1$-separable
and in standard form. Moreover, assuming PFA, filtration-equivalent $\aleph
_1$-separable groups of cardinality $\aleph _1$ are isomorphic. Thus the
result follows from Theorem \ref{stform}.

\section{A diamond construction: one tree}

The result to be proved in this section is the following:

\begin{theorem}
\label{onetree} Assume $\diamondsuit $. For any bounded tree $T_1$ there
exist non-isomorphic $\aleph _1$-separable groups $G^0$ and $G^1$ of
cardinality $\aleph _1$ which are $T_1$-equivalent (and
filtration-equivalent) and are both  in
standard form.
\end{theorem}

\TeXButton{Proof}{\proof} We will present the proof in layers of increasing
detail.

\smallskip\ 

{\bf (I)} Fix a stationary subset $E$ of $\omega _1$ consisting of limit
ordinals and such that $E$ is the disjoint union of two uncountable subsets $%
E_0$ and $E_1$ such that $\diamondsuit (E_1)$ holds.

Given a bounded tree $T$ (which in practice will be determined by, but not
equal to, $T_1)$, we shall identify its nodes with countable ordinals in
such a way that if $\nu <_T\mu $ (in the tree ordering), then $\nu <\mu $
(as ordinals).

By induction on $\alpha <\omega _1$ we will define the following data:

\begin{enumerate}
\item  continuous chains $\{G_\nu ^\ell :\nu <\alpha \}$ of countable free
groups (for $\ell =0,1$) such that for all $\nu <\mu <\alpha $, $G_\mu ^\ell
/G_\nu ^\ell $ is free if $\nu \notin E_1$, and if $\nu \in E_1$, then $%
G_{\nu +1}^\ell /G_\nu ^\ell $ has rank at most 1.

\item  homomorphisms $\pi _{\nu ,\mu }^\ell :G_\mu ^\ell \rightarrow G_\nu
^\ell $ for $\nu \leq \mu <\alpha $ and $\nu \notin E_1$ such that: $\pi
_{\nu ,\mu }^\ell $ is the identity on $G_\nu ^\ell $ ; for $\nu \leq \mu
<\rho $, $\pi _{\nu ,\mu }^\ell \subseteq \pi _{\nu ,\rho }^\ell $; and for $%
\tau <\nu \leq \mu $, $\pi _{\tau ,\nu }^\ell \circ \pi _{\nu ,\mu }^\ell
=\pi _{\tau ,\mu }^\ell $

\noindent (i.e., $\pi _{\nu ,\mu }^\ell $ is a projection and the system of
projections is coherent);

\item  for each $\nu $ with $\nu +1<\alpha $ an isomorphism $f_\nu ^0:G_{\nu
+1}^0\rightarrow G_{\nu +1}^1$ satisfying:

\begin{quote}
if $\nu _1<_T\nu _2$, then $f_{\nu _2}^0 \rest G_{\nu _1+1}^0=f_{\nu _1}^0$.
\end{quote}
\end{enumerate}

\noindent (These partial isomorphisms will give $\exists $ her winning
strategy.)

For convenience we will use $f_\nu ^1$ to denote $(f_\nu ^0)^{-1}:G_{\nu
+1}^1\rightarrow G_{\nu +1}^0$.

\medskip

Define $G^\ell =\cup _{\nu <\omega _1}G_\nu ^\ell $. (It depends on $T$, but
we suppress that in the notation.) Now we will indicate how we choose $T$ so
that $G^0$ and $G^1$ are $T_1$-equivalent.

\ Let $T_2={}^{<\omega _1}\omega _1\setminus \emptyset $, i.e., the tree of
non-empty countable sequences of countable ordinals, partially ordered by
inclusion (so it has $\aleph _1$ nodes of height $0$). Let $T$ be the
product $T_1\otimes T_2$, i.e., the (bounded) tree whose nodes are elements $%
(s,\sigma )\in T_1\times T_2$, where $s$ and $\sigma $ have the same height,
and the partial ordering is defined coordinate-wise. (As above, we identify
the nodes of $T$ with ordinals.)

Suppose we are able to carry out the construction outlined above for this $T$%
. Then since the $G_\nu ^\ell $ are free, $G^\ell $ is $\aleph _1$-free.
Moreover, for $\nu \notin E_1$ $\bigcup_{\mu <\omega _1}\pi _{\nu ,\mu
}^\ell :G^\ell \rightarrow G_\nu ^\ell $ is a projection which shows that $%
G_\nu ^\ell $ is a direct summand of $G^\ell $ ; so $G^\ell $ is $\aleph _1$%
-separable (and has a coherent system of projections; the fact that it is in
standard form will follow from the details of the construction --- see part
(V)).

We claim that $G^0$ and $G^1$ are $T_1$-equivalent. In fact, here is $%
\exists $'s winning strategy in the $T_1$-game. If in his first move $%
\forall $ plays $s_0\in T_1$ (which we may assume has height $0$), and $%
y_0\in G_{\gamma _0}^{\ell _0}$, $\exists $ chooses $\alpha _0$ such that $%
(s_0,\left\langle \alpha _0\right\rangle )\in T$ is the element $\nu _0$ in
the enumeration of $T$, where $\nu _0\geq \gamma _0$; and she plays $f_{\nu
_0}^{\ell _0}(y_0)\in G_{\nu _0+1}^{1-\ell _0}$. (Note that the domain of $%
f_{\nu _0}^{\ell _0}$ is $G_{\nu _0+1}^{\ell _0}\supseteq G_{\gamma
_0}^{\ell _0}$.) Suppose that after $\beta $ moves $\forall $ has chosen $%
s_0<_{T_1}s_1<_{T_1}...<_{T_1}s_\iota <_{T_1}...$ in the tree and $%
y_0,y_1,...,y_\iota ,...$ in the groups where $y_\iota \in G^{\ell _\iota 
\hbox{ }}$ ($\iota <\beta $), and $\exists $ has responded to the $\iota $th
move with $f_{\nu _\iota }^{\ell _\iota }(y_\iota )$ where $\nu _\iota
=(s_\iota ,\left\langle \alpha _0,...,\alpha _\iota \right\rangle )$. Now if 
$\forall $ plays $s_\beta >_{T_1}s_\iota $ ($\iota <\beta $) --- which we
can assume has height $\beta $ --- and $y_\beta \in G_{\gamma _\beta }^{\ell
_\beta }$, then $\exists $ chooses $\alpha _\beta $ such that $(s_\beta
,\left\langle \alpha _0,...,\alpha _\beta \right\rangle )$ is $\nu _\beta
\geq \gamma _\beta $, and plays $f_{\nu _\beta }^\ell (y_\beta )$. Notice
that $\nu _\beta >_T\nu _\iota $, so $f_{\nu _\beta }^\ell $ extends $f_{\nu
_\iota }^\ell $ for $\ell =0,1$. Therefore the sequence of moves determines
a partial isomorphism, so $\exists $ will win.

\medskip\ 

{\bf (II)} Of course, we also want to do the construction so that $G^0$ and $%
G^1$ are not isomorphic. This will be achieved by our construction of $%
G_{\delta +1}^\ell $ for $\delta \in E_1$ (plus the requirement 4 below);
when $\delta \in E_1$ we will make use of the ``guess'' provided by $%
\diamondsuit (E_1)$ of an isomorphism: $G_\delta ^0\rightarrow G_\delta ^1$.

Our construction will be such that when $\alpha =\mu +1$ where $\mu \notin E$%
, then 
$$
G_\alpha ^\ell =G_\mu ^\ell \oplus {\Bbb Z}x_{\mu ,0}^\ell \oplus {\Bbb Z}%
x_{\mu ,1}^\ell 
$$
\noindent When $\alpha =\sigma +1$ where $\sigma \in E_0$, then 
$$
G_\alpha ^\ell =G_\sigma ^\ell \oplus \bigoplus_{n\in \omega }{\Bbb Z}%
u_{\sigma ,n}^\ell \oplus {\Bbb Z}v_{\sigma ,n}^\ell \hbox{.} 
$$
\noindent We define 
$$
w_{\sigma ,n}=2u_{\sigma ,n+1}^0-u_{\sigma ,n}^0\hbox{. } 
$$

Notice that $\{w_{\sigma ,n}:n\in \omega \}$ generates a pure subgroup of $%
\bigoplus_{n\in \omega }{\Bbb Z}u_{\sigma ,n}^0$ which is not a direct
summand. Hence there is no isomorphism of $\bigoplus_{n\in \omega }{\Bbb Z}%
u_{\sigma ,n}^0\oplus {\Bbb Z}v_{\sigma ,n}^0$ with $\bigoplus_{n\in \omega }%
{\Bbb Z}u_{\sigma ,n}^1\oplus {\Bbb Z}v_{\sigma ,n}^1$ which takes each $%
w_{\sigma ,n}$ to $v_{\sigma ,n}^1$. In order to carry out the inductive
construction we will define in addition:

\begin{quote}
4. subsets $W_\alpha [\Theta ]$ of $G_\alpha ^0$ for every non-empty finite
subset $\Theta $ of $\alpha $ which is an antichain in $T$, satisfying:

(a) for all $\alpha <\beta $, $W_\alpha [\Theta ]\subseteq W_\beta [\Theta ]$%
;

(b) every element of $W_\alpha [\Theta ]$ is of the form $w_{\sigma ,n}$ for
some $\sigma \in E_0$, and $n\in \omega$. 
\end{quote}

\noindent The functions $f_\alpha ^0$ will be required to satisfy:

\begin{quote}
(c) for all $\mu \leq \alpha $, $j\in \{0,1\}$ $f_\alpha ^0(x_{\mu
,j}^0)=x_{\mu ,j}^1$; moreover, if $w_{\sigma ,n}\in W_{\alpha +1}[\Theta ]$
and $\Theta \cap \{\nu :\nu \leq _T\alpha \}\neq \emptyset $, then $f_\alpha
^0(w_{\sigma ,n})=v_{\sigma ,n}^1$.
\end{quote}

\noindent For any finite antichain $\Theta $ in $T$, let $W[\Theta
]=\bigcup_\alpha W_\alpha [\Theta ]$.

\medskip

Now we will outline how we do the construction so that $G^0$ and $G^1$ are
not isomorphic. Before we start, we choose a function $\Upsilon $ with
domain $E_0$ which maps onto the set of all $\omega $-sequences $%
\left\langle \Theta _n:n\in \omega \right\rangle $ of finite subsets of $T$
such that $\bigcup_{n\in \omega }\Theta _n$ is an antichain; we also require
that if $\Upsilon (\sigma )=\left\langle \Theta _n^\sigma :n\in \omega
\right\rangle $, then each $\Theta _n^\sigma \subseteq \sigma $.

Suppose now that we have defined $G_\nu ^\ell $ for $\nu \leq \alpha $. If $%
\alpha =\sigma \in E_0$, then $G_{\sigma +1}^\ell $ will be defined as indicated
above and is such that (as we will prove)

\begin{quote}
(II.1) for all $e\in \{1,-1\}$, there is no
isomorphism of $G_{\sigma +1}^0$ with $\bigoplus_{n\in \omega }{\Bbb Z}%
v_{\sigma ,n}^1\oplus C $ for any $C$, which for all $n\in \omega $ takes $%
w_{\sigma ,n}$ to $ev_{\sigma ,n}^1$.
\end{quote}

\noindent  Moreover $w_{\sigma ,n}$ will be put into $W_{\sigma +1}[\Theta
_n^\sigma ]$. (This is the only way that an element becomes a member of a $%
W_\alpha [\Theta ]$.)

If $\alpha =\delta \in E_1$ and $\beta <\delta $, we introduce the notation $%
A_{\beta ,\delta }=\{t:t$ is $<_T$-minimal in $\delta \setminus \beta \}$
--- so $A_{\beta ,\delta }$ is an antichain. We fix finite subsets $\Theta
_n^{\beta ,\delta }$ of $A_{\beta ,\delta }$ which form a chain such that $%
\cup _{n\in \omega }\Theta _n^{\beta ,\delta }=A_{\beta ,\delta }$. We
consider the prediction given by $\diamondsuit (E_1)$ of an isomorphism $h:$ 
$G_\delta ^0\rightarrow G_\delta ^1$ and we ask whether the following holds:

\begin{quote}
(II.2) $\exists $ $\beta <\delta $ $\forall $ $e\in \{1,-1\}$ $\forall $ $%
n\in \omega $ $\exists $ $w_{\sigma ,m}\in W_\delta [\Theta _n^{\beta
,\delta }]$ such that $h(w_{\sigma ,m})\neq ev_{\sigma ,m}^1$.
\end{quote}

We will do the construction of $G_{\delta +1}^\ell $ so that:

\begin{quote}
(II.3) If (II.2) holds, then $G_{\delta +1}^\ell
/G_\delta ^\ell $ is non-free rank 1 and  $h$ does not extend to a
homomorphism: $G_{\delta +1}^0\rightarrow G_{\delta +1}^1$.
\end{quote}

\medskip 

Assuming we can do all of this, let us see why $G^0$ is not isomorphic to $%
G^1$. Suppose, to the contrary, that there is an isomorphism $%
H:G^0\rightarrow G^1$. Then there is a stationary set, $S$, of $\delta \in
E_1$ where $\diamondsuit (E_1)$ guesses $h=H \rest G_\delta ^0$ and $H:G_\delta
^0\rightarrow G_\delta ^1$.  Note that Lemma \ref{stein} implies that $H$
must map $G_{\delta +1}^0$ into $G_{\delta +1}^1$ because $G_{\delta
+1}^0/G_\delta ^0$ is non-free rank $1$ but $G^1/G_{\delta +1}^1$ is $\aleph
_1$-free by construction.  If for any such $\delta $ (II.2) holds, then $%
H \rest G_{\delta +1}^0$ would extend $h=H \rest G_\delta ^0$, contradicting (II.3).

Since (II.2) fails, for all $\delta \in S$ and all $\beta <\delta $ there
exists $e\in \{1,-1\}$ and a finite subset $\Theta $ of $A_{\beta ,\delta }$
such that $H(w_{\sigma ,n})=$ $ev_{\sigma ,n}^1$ for all $w_{\sigma ,n}\in
W_\delta [\Theta ]$. Now there is a cub $C$ such that for all $\delta \in C$%
, all $e\in \{1,-1\}$, all $\beta <\delta $, and all finite subsets $\Theta $
of $A_{\beta ,\delta }$, if $H(w_{\sigma ,n})\neq $ $ev_{\sigma ,n}^1$ for
some $w_{\sigma ,n}\in W[\Theta ]$, then $H(w_{\sigma ,n})\neq $ $ev_{\sigma
,n}^1$ for some $w_{\sigma ,n}\in W_\delta [\Theta ]$. Thus for all $\delta
\in C\cap S$ and all $\beta <\delta $, there exists $e\in \{1,-1\}$ and a
finite subset $\Theta $ of $A_{\beta ,\delta }$ such that $H(w_{\sigma ,n})=$
$ev_{\sigma ,n}^1$ for all $w_{\sigma ,n}\in W[\Theta ]$. Since $C\cap S$ is
uncountable, it follows easily that there exists $e\in \{1,-1\}$, and an
uncountable set $\{\Theta _\nu :\nu <\omega _1\}$ of pairwise disjoint
finite antichains such that $H(w_{\sigma ,n})=$ $ev_{\sigma ,n}^1$ for all $%
w_{\sigma ,n}\in W[\Theta _\nu ]$ for all $\nu <\omega _1$. Since $T$ has no
uncountable branches, by a standard argument (see, for example, 
\cite[Lemma 24.2, p. 245]{J}), there is a countably infinite subset $\{\nu
_n:n\in \omega \}$ of $\omega _1$ such that $\bigcup \{\Theta _{\nu _n}:n\in
\omega \}$ is an antichain. There exists $\sigma \in E_0$ such that $%
\Upsilon (\sigma )=\left\langle \Theta _{\nu _n}:n\in \omega \right\rangle $%
. Now $H \rest G_{\sigma +1}^0$ is such that for all $n\in \omega $, $H(w_{\sigma
,n})=$ $ev_{\sigma ,n}^1$ since $w_{\sigma ,n}\in W_{\sigma +1}[\Theta _{\nu
_n}]$; this contradicts (II.1), since $\bigoplus_{n\in \omega }{\Bbb Z}%
v_{\sigma ,n}^1$ is a direct summand of $G_{\sigma +1}^1$, and hence of $G^1$
(by 2).
\medskip\ 

{\bf (III)} The next step is to describe in detail the recursive
construction of the data satisfying the properties 1, 2, 3 and 4, as well as
(II.1) and (II.3). So assume that we have defined $G_\nu ^\ell $%
, and $W_\nu [\Theta ]$ for $\nu $ $<\alpha $ and $f_\nu ^\ell $ for $\nu
+1<\alpha $.

There are several cases to consider.

\smallskip\ 

{\it Case 1: }$\alpha ${\it \ is a limit ordinal.} We let $G_\alpha ^\ell
=\cup _{\nu <\alpha }G_\nu ^\ell $, $W_a[\Theta ]=\bigcup_{\nu <\alpha
}W_\nu [\Theta ]$. Clearly the desired properties are satisfied.

\smallskip\ 

If $\alpha $ is a successor, $\alpha =\mu +1$, we will define $G_\alpha
^\ell $ so that

\begin{quote}
(III.1) if $B=\{t:t<_T\mu \}$ and we define $g_B=\cup \{f_t^0:t\in B\}$,
then $g_B$ (which is a function by 3.) extends to an isomorphism, $f_\mu ^0$%
, of $G_\alpha ^0$ onto $G_\alpha ^1$ which satisfies 4(c), i.e.  for all $%
\nu \leq \mu $, $j\in \{0,1\}$ $f_\mu ^0(x_{\nu ,j}^0)=x_{\nu ,j}^1$ and  if 
$w_{\sigma ,n}\in W_\alpha [\Theta ]$ and $\Theta \cap \{\nu :\nu \leq _T\mu
\}\neq \emptyset $, then $f_\mu ^0(w_{\sigma ,n})=v_{\sigma ,n}^1$.
\end{quote}

\ 

Leaving the verification of (III.1) to the next part, we will show how to
define the data at $\alpha $ (except for the definition of the $\pi _{\sigma
,\alpha }^\ell $ which we defer to part (V)).

\smallskip\ 

{\it Case 2: }$\alpha =\mu +1$ {\it for some }$\mu \notin E.$ As described
above, define 
$$
G_\alpha ^\ell =G_\mu ^\ell \oplus {\Bbb Z}x_{\mu ,0}^\ell \oplus {\Bbb Z}%
x_{\mu ,1}^\ell \hbox{.} 
$$
Let $W_\alpha [\Theta ]=W_\mu [\Theta ]$ for every $\Theta \subseteq \mu $ ($%
=\emptyset $ if $\Theta $ is not a subset of $\mu $). Assuming (III.1), we
have $f_\mu ^0$ as desired$.$

\smallskip\ 

{\it Case 3: }$\alpha =\sigma +1${\it , where }$\sigma \in E_0${\it . }In
this case, as stated before, 
$$
G_\alpha ^\ell =G_\sigma ^\ell \oplus \bigoplus_{n\in \omega }{\Bbb Z}%
u_{\sigma ,n}^\ell \oplus {\Bbb Z}v_{\sigma ,n}^\ell 
$$
and recall that $w_{\sigma ,n}$ is defined to be $2u_{\sigma
,n+1}^0-u_{\sigma ,n}^0$. Say $\Upsilon (\sigma )=\left\langle \Theta
_n^\sigma :n\in \omega \right\rangle $. Define 
$$
W_\alpha ^{}[\Theta ]=\left\{ 
\begin{array}{ll}
W_\sigma [\Theta ]\cup \{w_{\sigma ,n}\} & \hbox{if }\Theta =\Theta
_n^\sigma \\ W_\sigma [\Theta ] & \hbox{otherwise.} 
\end{array}
\right. 
$$

Assuming (III.1) (with $\mu =\sigma $), we can define $f_\sigma ^0$. Now let
us see why (II.1) holds. Suppose to the contrary that there is an
isomorphism $H:G^0\rightarrow G^1$ contradicting (II.1). Now $%
\bigoplus_{n\in \omega }{\Bbb Z}v_{\sigma ,n}^1$ is a direct summand of $%
G_\alpha ^1$ and hence (by 2) a direct summand of $G^1$. Thus $%
H^{-1}[\bigoplus_{n\in \omega }{\Bbb Z}v_{\sigma ,n}^1]$ is a direct summand
of $G^0$. But by assumption on $H$, $H^{-1}[\bigoplus_{n\in \omega }{\Bbb Z}%
v_{\sigma ,n}^1]=\bigoplus_{n\in \omega }{\Bbb Z}w_{\sigma ,n}^{}$ and the
latter is {\it not} a direct summand of $G^0$ because the coset of $%
u_{\sigma ,0}^0$ is a non-zero element of $G^0/\bigoplus_{n\in \omega 
}{\Bbb Z}
w_{\sigma ,n}^{}$ which is divisible by all power of $2$ by definition of
the $w_{\sigma ,n}$.

\smallskip\ 

{\it Case 4: }$\alpha =\delta +1${\it , where }$\delta \in E_1$. If (II.2)
fails, let $G_{\delta +1}^\ell =G_\delta ^\ell $. Otherwise, let $\beta $ be
as in (II.2). We introduce some {\it ad hoc }notation. For any finite subset 
$\Theta $ of $A_{\beta ,\delta }$, let $f_\Theta $ be the function whose
domain is the subgroup generated by $\{x_{\mu ,j}^0:\mu \notin E$, $\mu
<\delta $, $j\in \{0,1\}\}\cup W_\delta [\Theta ]$ such that $f_\Theta
(x_{\mu ,j}^0)=x_{\mu ,j}^1$ and $f_\Theta (w_{\sigma ,n})=v_{\sigma ,n}^1$.
Notice that for all $u\in \limfunc{dom}(f_\Theta )$ and all $\nu \in \Theta $%
, if $\nu \leq _T\rho $ and $u\in \limfunc{dom}(f_\rho ^0)$, then $f_\Theta
(u)=f_\rho ^0(u)$ by 4(c). Let $\Theta _n^{\beta ,\delta }$ be as before
(finite subsets forming a chain whose union is $A_{\beta ,\delta }$); for
short, let $\Theta _n=\Theta _n^{\beta ,\delta }$. We claim that:

\begin{quote}
(III.2) given $m,m^{\prime }\in {\Bbb Z} \setminus \{0\}$, $n\in \omega $, $%
y\in G_\delta ^1$, for sufficiently large $\gamma $ $<\delta $ there exists $%
k_{}^0\in \limfunc{dom}(f_{\Theta _n})\cap G_{\gamma +2}^0$ such that $k^0$
is pure-independent mod $G_{\gamma +1}^0$ and is such that $mh(k^0)\neq
m^{\prime }f_{\Theta _n}(k^0)+y$. Moreover, $f_{\Theta _n}(k^0)$ is
pure-independent mod $G_{\gamma +1}^1$.
\end{quote}

Supposing this is true --- we will prove it in part (IV) --- let us define $%
G_{\delta +1}^\ell $. Fix a ladder $\eta _\delta $ on $\delta $. Also,
enumerate in an $\omega $-sequence all triples $\left\langle
r,d,v\right\rangle $ where $r\in \omega $, $d\in {\Bbb Z} \setminus \{0\}$,
and $g\in G_\delta ^1$ so that the $n$th triple $\left\langle
r,d,g\right\rangle $ satisfies $n>r$. By (III.2) we can inductively define
primes $p_n$, ordinals $\gamma _n\geq \eta _\delta (n)$, and elements $%
k_{\delta ,n}^0\in \limfunc{dom}(f_{\Theta _n})\cap G_{\gamma _n+2}^0$
pure-independent over $G_{\gamma _n+1}^0$ such that (if the $n$th triple is $%
\left\langle r,d,g\right\rangle $), $p_n$ does not divide $mh(k_{\delta
,n}^0)-m^{\prime }f_{\Theta _n}(k_{\delta ,n}^0)-y$ where 
$$
\begin{array}{c}
m=\prod_{i=0}^{n-1}p_i \\ 
m^{\prime }=d\prod_{i=r}^{n-1}p_i \\ 
y=\sum_{j=0}^{n}(\prod_{i=0}^{j-1}p_i)h(k_{\delta
,j}^0)+g-d\sum_{j=r}^{n}(\prod_{i=r}^{j-1}p_i)f_{\Theta _j}(k_{\delta
,j}^0)\hbox{.}
\end{array}
$$
(Note that since $G_\delta ^1$ is free, every  non-zero element is divisible
by only finitely many primes, so we can take $p_n$ to be any sufficiently
large prime.) Then we let $G_{\delta +1}^0$ be generated by $G_\delta ^0\cup
\{z_{\delta ,n}^0:n\in \omega \}$ modulo the relations%
$$
p_nz_{\delta ,n+1}^0=z_{\delta ,n}^0+k_{\delta ,n}^0 
$$
and $G_{\delta +1}^1$ is defined similarly, except that we impose the
relations 
$$
p_nz_{\delta ,n+1}^1=z_{\delta ,n}^1+f_{\Theta _n}(k_{\delta ,n}^0)\hbox{.} 
$$

We need to show that $h$ does not extend to a homomorphism: $G_{\delta
+1}^0\rightarrow G_{\delta +1}^1$. If it does, then $h(z_{\delta
,0}^0)=dz_{\delta ,r}^1+g$ for some $r\in \omega $, $d\in {\Bbb Z} \setminus
\{0\}$, and $g\in G_\delta ^1$. Let $n$ be such that $\left\langle
r,d,g\right\rangle $ is the $n$th triple in the list. Now, in $G_{\delta
+1}^0$ we have%
$$
(\prod_{i=0}^np_i)z_{\delta ,n+1}^0=z_{\delta
,0}^0+\sum_{j=0}^n(\prod_{i=0}^{j-1}p_i)k_{\delta ,j}^0 
$$
so, applying $h$, we conclude that $p_n$ divides 
$$
dz_{\delta ,r}^1+g+\sum_{j=0}^n(\prod_{i=0}^{j-1}p_i)h(k_{\delta ,j}^0)\hbox{%
.} 
$$
On the other hand, in $G_{\delta +1}^1$ we have $p_n$ divides%
$$
dz_{\delta ,r}^1+d\sum_{j=r}^n(\prod_{i=r}^{j-1}p_i)f_{\Theta _j}(k_{\delta
,j}^0) 
$$
so, subtracting, we obtain a contradiction since $p_n$ divides $mh(k_{\delta
,n}^0)-m^{\prime }f_{\Theta _n^\sigma }^0(k_{\delta ,n}^0)-y$, where $m$, $%
m^{\prime }$, and $y$ are as above.

We let $W_{\delta +1}[\Theta ]=W_\delta [\Theta ]$ for any subset $\Theta $
of $\delta $ (and $=\emptyset $ if $\Theta \not \subseteq \delta $). By
(III.1) we can define $f_\delta ^0$.

\smallskip\ 

\medskip\ 

{\bf (IV)} In this layer we will prove (III.1) and (III.2).

First let us prove (III.2) since for the purposes of proving (III.1) we will
need more information about the nature of the elements $k_{\delta ,n}^0$.
Fix $m,m^{\prime },n,y,\gamma $ as in (III.2); there are several cases. In
the first two cases we can use any $\gamma <\delta $.

{\it Case (i): }$y\neq 0${\it .} If neither $x_{\gamma +1,0}^0$ nor $%
x_{\gamma +1,1}^0$ will serve for $k^0$, then $x_{\gamma +1,0}^0-x_{\gamma
+1,1}^0$ will.

\smallskip\ 

{\it Case (ii): }$y=0$, $m\neq \pm m^{\prime }${\it . }Let{\it \ }$%
k^0=x_{\gamma +1,0}^0$. Then by construction, $k^0$ generates a cyclic
summand of $G_\delta ^0$; hence $f_{\Theta _n}(k^0)$ and $h(k^0)$ both
generate cyclic summands of $G_\delta ^1$. Hence  $mh(k^0)\neq m^{\prime
}f_{\Theta _n}^0(k^0)$.

\smallskip\ 

{\it Case (iii): }$y=0${\it , }$m=m^{\prime }${\it .} Pick $\gamma $
sufficiently large so that there exists $w_{\sigma ,j}\in G_{\gamma +1}^0$ $%
\cap W_\delta [\Theta _n]$ such that $f_{\Theta _n}(w_{\sigma ,j})\neq
h(w_{\sigma ,j})$. If $x_{\gamma +1,0}^0$ will not serve for $k^0$ (i.e., $%
h(x_{\gamma +1,0}^0)=x_{\gamma +1,0}^1$), then we can take $k^0$ to be $%
x_{\gamma +1,0}^0+w_{\sigma ,j}$.

\smallskip\ 

{\it Case (iv): }$y=0${\it , }$m=-m^{\prime }$. Similarly $k^0$ can be taken
to be of the form $x_{\gamma +1,0}^0$ or $x_{\gamma +1,0}^0- w_{\sigma ,j} $ where $%
f_{\Theta _n}(w_{\sigma ,j} )\neq -h(w_{\sigma ,j} )$.

\smallskip\ 

Now if we examine the construction in Case 4 of (III) and the proof above we
see that

\begin{quote}
(IV.1) each $k_{\delta ,n}^0$ can be (and will be) taken to be of the form $%
x_{\mu _n,j_n}^0\pm \xi _{\delta ,n}$ where $\xi _{\delta ,n}$ is $0$, $%
x_{\sigma ,j}^0$ or $w_{\sigma ,j}$ for some $\sigma ,j$.
\end{quote}

\noindent We will say that $w_{\sigma ,j}$ {\it is a part of} $k_{\delta
,n}^0$ in case $\xi _{\delta ,n}$ is $w_{\sigma ,j}$.

\smallskip

Before beginning the proof of (III.1), let us observe the following facts:

\begin{quote}
(IV.2) Given $\sigma \in E_0$ and $N\in \omega $, there is an isomorphism $%
g^{\prime }:\bigoplus_{n\in \omega }{\Bbb Z}u_{\sigma ,n}^0\oplus {\Bbb Z}%
v_{\sigma ,n}^0\rightarrow \bigoplus_{n\in \omega }{\Bbb Z}u_{\sigma
,n}^1\oplus {\Bbb Z}v_{\sigma ,n}^1$ such that $g^{\prime }(w_{\sigma
,n})=v_{\sigma ,n}^1$ for $n\leq N$ and $g^{\prime }(u_{\sigma
,n}^0)=u_{\sigma ,n}^1$ for $n\geq N+1$.
\end{quote}

\noindent Indeed, we can define $g^{\prime }(u_{\sigma ,n}^0)=2g^{\prime
}(u_{\sigma ,n+1}^0)-v_{\sigma ,n}^1$ for $n\leq N$ (and the other values
appropriately).

\begin{quote}
%Given an isomorphism $g:G_\sigma ^0\rightarrow G_\sigma ^1$ where $%
%\sigma \in E_0$, and given $N\in \omega $, we can extend $g$ to an
%isomorphism $g^{\prime }:$ $G_{\sigma +1}^0\rightarrow G_{\sigma +1}^1$ such
%that $g^{\prime }(w_{\sigma ,n})=v_{\sigma ,n}^1$ for $n\leq N$. In fact, we
%can define $g^{\prime }(z_{\sigma ,n}^0)=z_{\sigma ,n}^1$ for $n\geq N+1$,
%and $g^{\prime }(z_{\sigma ,n}^0)=2g^{\prime }(z_{\sigma ,n+1}^0)-v_{\sigma
%,n}^1$ for $n\leq N$ (and the other values appropriately).

(IV.3) Given an isomorphism $g:G_\delta ^0\rightarrow G_\delta ^1$ where $%
\delta \in E_1$, we can extend $g$ to an isomorphism $g^{\prime }:$ $%
G_{\delta +1}^0\rightarrow G_{\delta +1}^1$ provided that (using the
notation of Case 4) $g(k_{\delta ,n}^0)=f_{\Theta _n}(k_{\delta ,n}^0)$ for
almost all $n\in \omega $.
\end{quote}

\noindent Indeed, if $g(k_{\delta ,n}^0)=f_{\Theta _n}(k_{\delta ,n}^0)$ for
all $n\geq N$, we can define $g^{\prime }(z_{\delta ,n}^0)=z_{\delta ,n}^1$
for $n\geq N$ and $g^{\prime }(z_{\delta ,n}^0)=p_ng^{\prime }(z_{\delta
,n+1}^0)-g(k_{\delta ,n}^0)$ for $n<N$ by ``downward induction''. We will
apply (IV.3) to the situation of (III.1), with $g=g_B$, $\delta =\mu $, $%
\delta +1=\alpha $; if we are in Case 4, then the hypothesis on $g$ in
(IV.3) will hold if there exists $t\in B$ such that $t\geq \beta $ (where $%
\beta $ is as in Case 4).

We return to the notation of (III.1). Let $\tau =\sup \{t+1:t\in B\}$; then $%
\limfunc{dom}g_B=G_\tau ^0$. Assume first that $\tau =\mu $. In case $G_{\mu
+1}^\ell /G_\mu ^\ell $ is free there is no problem extending $g_B$; in the
other case $\mu =\delta \in E_1$ and by the remarks above we can extend $g_B$
since there exists $t\in B$ such that $t\geq \beta $ (since $\sup B=\delta $%
).

We are left with the case when $\tau <\mu $. We will first define an
extension of $g_B$ to a partial isomorphism $\tilde g_B$ whose domain is 
$$
\limfunc{dom}(g_B)+\left\langle 
\begin{array}{c}
\{x_{\nu ,j}^0:\nu 
\hbox{ }\leq \mu \hbox{, }j=0,1\}\cup \\ \{u_{\sigma ,n}^0:\sigma \in
E_0\cap \mu +1 
\hbox{, }n\in \omega \}\cup \\ \{v_{\sigma ,n}^0:\sigma \in E_0\cap \mu +1%
\hbox{, }n\in \omega \} 
\end{array}
\right\rangle 
$$

Notice that every $k_{\delta ,n}^0$ for $\delta \leq \mu ,n\in \omega $
belongs to the domain of $\tilde g_B$. We let $\tilde g_B(x_{\nu
,j}^0)=x_{\nu ,j}^1$ for all $\nu ,j$. By enumerating in an $\omega $%
-sequence the set $(E_0\cup E_1)\cap (\mu +1)$ we can define by recursion the
values $\tilde g_B(u_{\sigma ,n}^0)$ and $\tilde g_B(v_{\sigma ,n}^0)$  so that:

\begin{itemize}
\item  $\tilde g_B(w_{\sigma ,n})=v_{\sigma ,n}^1$ whenever $w_{\sigma
,n}\in W_{\mu +1}[\Theta ]$ for some $\Theta $ with $B\cap \Theta \neq
\emptyset $;

\item  for all $\sigma \in E_0$ with $\tau \leq \sigma \leq \mu $, for
almost all $n\in \omega $, $\tilde g_B(u_{\sigma ,n}^0)=u_{\sigma ,n}^1$; and

\item  for all $\delta \in E_1$ with $\tau \leq \delta \leq \mu $, for
almost all $n\in \omega $, if (for some $\sigma ,m$) $w_{\sigma ,m}$ is a
part of $k_{\delta ,n}^0$ , then $\tilde g_B(w_{\sigma ,m})=v_{\sigma ,m}^1$.
\end{itemize}

\noindent The first condition is required by 4(c). In view of (IV.2), there
is no conflict between the first two conditions because for any $\sigma \in
E_0$, $\bigcup_{n\in \omega }\Theta _n^\sigma $ is an antichain, so there is
at most one $n$ such that $\Theta _n^\sigma \cap B\neq \emptyset $.

To be sure that the third condition can indeed be satisfied, we need to
consider the case that for some $\delta \in E_1$, there are infinitely many $%
n$ such that there exists $w_{\sigma _n,m_n}$ which is a part of $k_{\delta
,n}^0$ and belongs to the domain of $g_B$. Say this is the case for $n$
belonging to the (infinite) set $Y\subseteq \omega $ (for a fixed $\delta $%
). Then for each $n\in Y$ $\exists t_n\in B$ such that $t_n\geq \sigma _n$.
Suppose that the construction of $G_{\delta +1}^\ell $ uses $A_{\beta
,\delta }=\cup _{n\in \omega }\Theta _n^{\beta ,\delta }$. Selecting one $%
n_{*}\in Y$, we see that since $\Theta _{n_{*}}^{\beta ,\delta }\subseteq
\sigma _{n_{*}}$, $\sigma _{n_{*}}>\beta $ and hence $t_{n_{*}}\in A_{\beta
,\delta }$. Therefore there exists $M$ such that for all $n\geq M$, $%
t_{n_{*}}\in \Theta _n^{\beta ,\delta }$. But then, for $n\in Y$ with $n\geq
M$, $t_n\geq \sigma _n\supseteq \Theta _n^{\beta ,\delta }$, so $%
t_{n_{*}}\leq t_n$ and thus $t_{n_{*}}\leq _Tt_n$. By the construction in
Case 4 and by 4(c), $g_B(w_{\sigma _n,m_n})=v_{\sigma _n,m_n}^1$ for $n \in Y
$, $n\geq M$.
Moreover, there is no conflict between the last two conditions because, by
construction, if $\delta \in E_1$ and $\sigma \in E_0$, then $w_{\sigma
,m}\in W_\delta [\Theta _n^{\beta ,\delta }]$ if and only if $\Theta
_n^{\beta ,\delta }=\Theta _m^\sigma $, but the elements of $\{\Theta
_m^\sigma :m\in \omega \}$ are disjoint and the $\Theta _n^{\beta ,\delta }$
form a chain under $\subseteq .$

It remains to extend $\tilde g_B$ to $f^0_\mu$ by defining 
$f^0_\mu(z_{\delta ,n}^0)$ for $\tau \leq \delta \leq \mu ,n\in \omega $. This is
possible by observation (IV.3) because of the construction of $\tilde g_B$.

\medskip\ 

{\bf (V)} We will define the projections $\pi _{\nu ,\mu }^\ell $ by
induction on $\mu $  and then verify the conditions to be in
standard form (see section 1 or  \cite[Def. 1.9(ii), p. 257]{EM}). We refer to the cases
of the construction in part (III). In Case 1, we take unions. In Case 2, for 
$\nu <\mu +1$ we let $\pi _{\nu ,\mu +1}^\ell $ be the extension of $\pi
_{\nu ,\mu }^\ell $ which sends each $x_{\mu ,j}^\ell $ to 0 $.$ (Here, $\pi
_{\mu ,\mu }^\ell $ is the identity.) In Case 3, for $\nu \leq \sigma $ we
let $\pi _{\nu ,\sigma +1}^\ell $ be the extension of $\pi _{\nu ,\sigma
}^\ell $ which sends each $u_{\sigma ,n}^\ell $ and each $v_{\sigma ,n}^\ell 
$ to 0. 

Finally, for Case 4, we use the notation of that case. We define $\pi _{\nu
,\delta +1}^0(z_{\delta ,n}^0)=-\sum_{j=n}^md_{n,j}k_{\delta ,j}^0 $ where $%
m $ is maximal such that $\gamma _m+2\leq \nu $ and $d_{n,j}=%
\prod_{i=n}^{j-1}p_i $ (and $d_{n,0}=1 $ ). 
(Compare \cite[pp. 249f]{EM}.) The definition of $\pi _{\nu
,\delta +1}^1 $ is similar, replacing $k_{\delta ,j}^0 $ by $f_{\Theta
_j}(k_{\delta ,j}^0) $. Let $Y_\delta ^\ell =\{z_{\delta ,n}^\ell :n\in
\omega \} $. Then we can easily verify the conditions of 
\cite[Def. 1.9(ii), p. 257]{EM} using the information in the proof of
(III.2) about the form of $k_{\delta ,j}^0 $.

This completes the proof of Theorem 7.

\section{A non-structure theorem}

Our goal is to generalize the construction in the previous section to prove:

\begin{theorem}
\label{nonstr} Assume $\diamondsuit $. There exists an $\aleph _1$-separable
group $G^0$ and for each bounded tree $T_1$ an $\aleph _1$-separable group $%
G^{T_1}$ which is $T_1$-equivalent to $G^0$ but not isomorphic to $G^0$.
Moreover, all the groups are of cardinality $\aleph _1$ and in standard form.
\end{theorem}

\TeXButton{Proof}{\proof} We assume familiarity with the previous proof and
outline the modifications, in layers of increasing detail.

\smallskip\ 

{\bf (VI)} Fix a stationary subset $E$ of $\omega _1$ consisting of limit
ordinals ($>0$) and such that $E$ is the disjoint union of two subsets $E_0$
and $E_1$ such that cardinality $\diamondsuit (E_0)$ and $\diamondsuit (E_1) 
$ hold. ($\diamondsuit (E_0)$ is not essential, but convenient.)

We need only consider bounded trees $T$ on $\omega _1$ such that if $\nu
<_T\mu $ (in the tree ordering), then $\nu <\mu $ (as ordinals). For each $%
\delta \in E_1$ (resp. $\sigma \in E_0)$, diamond will give us a
``prediction'' $T_\delta =\left\langle \delta ,<_\delta \right\rangle $
(resp. $T_\sigma $) of the restriction of a bounded tree to $\delta $ (resp. 
$\sigma $). If $\mu <\delta $ we write $T_\delta  \rest \mu $ for $\left\langle
\mu ,<_\delta \cap (\mu \times \mu )\right\rangle $.

By induction on $\delta \in \{0\}\cup E$ we will define the following data:

\begin{enumerate}
\item  continuous chains $\{G_\nu ^\delta :\nu \leq \delta +1\}$ of
countable free groups such that for all $\nu <\mu \leq \delta +1$, $G_\mu
^\delta /G_\nu ^\delta $ is free if $\nu \notin E_1$, and if $\nu \in E_1$,
then $G_{\nu +1}^\delta /G_\nu ^\delta $ has rank at most 1.

\item  projections $\pi _{\nu ,\mu }^\delta :G_\mu ^\delta \rightarrow G_\nu
^\delta $ for $\nu \leq \mu \leq \delta +1$ and $\nu \notin E_1$ such that:
for $\nu \leq \mu <\rho $, $\pi _{\nu ,\mu }^\delta \subseteq \pi _{\nu
,\rho }^\delta $; and for $\tau <\nu \leq \mu $, $\pi _{\tau ,\nu }^\delta
\circ \pi _{\nu ,\mu }^\delta =\pi _{\tau ,\mu }^\delta $;

\item  for each $\delta \in E$ and each $\nu $ $\leq \delta $ an isomorphism 
$f_\nu ^\delta :G_{\nu +1}^0\rightarrow G_{\nu +1}^\delta $ satisfying:

\begin{quote}
if $\nu _1<_\delta \nu _2$, then $f_{\nu _2}^\delta  \rest G_{\nu _1+1}^0=f_{\nu
_1}^\delta $.
\end{quote}
\end{enumerate}

\smallskip\ \ 

Moreover, we require that if $\delta <\delta ^{\prime }$ are elements of $E$
such that $T_\delta =T_{\delta ^{\prime }} \rest \delta $, then $G_\nu ^\delta
=G_\nu ^{\delta ^{\prime }}$ for $\nu \leq \delta +1$; $\pi _{\nu ,\mu
}^{\delta ^{\prime }}=\pi _{\nu ,\mu }^\delta $ for $\nu \leq \mu \leq
\delta +1$; and $f_\nu ^{\delta ^{\prime }}=f_\nu ^\delta $ for $\nu \leq
\delta $.

\medskip

Define $G^0=\cup _{\nu <\omega _1}G_\nu ^0$ and for each bounded tree $T$ on 
$\omega _1$ let $G^T=\bigcup \{G_\nu ^\delta :T_\delta =T \rest \delta $, $\nu
\leq \delta +1\}$. As before, given $T_1$ we can choose $T$ so that $G^0$
and $G^T$ are $T_1$-equivalent.

We indicate how to modify the previous construction so that $G^0$ and $G^T$
are not isomorphic. Our construction will be such that when $\alpha =\mu +1$
where $\mu \notin E$, then 
$$
\hbox{(*) }G_\alpha ^0=G_\mu ^0\oplus {\Bbb Z}x_{\mu ,0}^0\oplus {\Bbb Z}%
x_{\mu ,1}^0 
$$
and 
$$
\hbox{(**) }G_\alpha ^\delta =G_\mu ^\delta \oplus {\Bbb Z}x_{\mu ,0}^1\oplus 
{\Bbb Z}x_{\mu ,1}^1 
$$
for $\delta \in E$, $\alpha <\delta $.

\noindent When $\alpha =\sigma +1$ where $\sigma \in E_0$, then 
$$
\hbox{(***) }G_\alpha ^0=G_\sigma ^0\oplus \bigoplus_{n\in \omega }{\Bbb Z}%
u_{\sigma ,n}^0\oplus {\Bbb Z}v_{\sigma ,n}^0 
$$
and 
$$
\hbox{(****) }G_\alpha ^\delta =G_\sigma ^\delta \oplus \bigoplus_{n\in
\omega }{\Bbb Z}u_{\sigma ,n}^1\oplus {\Bbb Z}v_{\sigma ,n}^1 
$$
for $\delta \in E$, $\sigma <\delta $.

\noindent We define 
$$
w_{\sigma ,n}=2u_{\sigma ,n+1}^0-u_{\sigma ,n}^0\hbox{. } 
$$

In order to carry out the inductive construction we will define in addition:

\begin{quote}
4. for $\delta \in E$ and $\alpha \leq \delta +1$, subsets $W_\alpha ^\delta
[\Theta ]$ of $G_\alpha ^0$ for every non-empty finite subset $\Theta $ of $%
\alpha $ which is an antichain in $T_\delta $, satisfying:

(a) for all $\alpha <\beta $, $W_\alpha ^\delta [\Theta ]\subseteq W_\beta
^\delta [\Theta ]$;

(b) every element of $W_\alpha ^\delta [\Theta ]$ is of the form $w_{\sigma
,n}$ for some $n\in \omega $ and some $\sigma \in E_0$   such that $T_\delta
 \rest \sigma =T_\sigma .$
\end{quote}

\noindent The functions $f_\alpha ^\delta $ will be required to satisfy (as
before):

\begin{quote}
(c) for all $\mu \leq \alpha $, $j\in \{0,1\}$ $f_\alpha ^\delta (x_{\mu
,j}^0)=x_{\mu ,j}^1$; moreover, if $w_{\sigma ,n}\in W_{\alpha +1}^\delta
[\Theta ]$ and $\Theta \cap \{\nu :\nu \leq _\delta \alpha \}\neq \emptyset $%
, then $f_\alpha ^\delta (w_{\sigma ,n})=v_{\sigma ,n}^1$.
\end{quote}

\noindent Moreover, in order to carry out the inductive construction we will
also require the following for all $\delta \in E$, $\alpha \leq \delta $ :%

\begin{quote}
(d) if $\sigma \in E_0$ with $\sigma \leq \alpha +1$ and $T_\delta  \rest \sigma
\neq T_\sigma $, then $f_\alpha ^\delta (u_{\sigma ,n}^0)=u_{\sigma ,n}^1$
for all $n\in \omega $;

(e) for all pairs $\beta _1$, $\beta _2$ with $\sup \{t:t<_\delta \alpha
\}\leq \beta _1<\beta _2\leq \alpha $, it is the case for almost all $n\in
\omega $ that for all $w_{\sigma ,m}\in W_{\alpha +1}^\delta [\Theta
_n^{\beta _1,\beta _2}]$ we have $f_\alpha ^\delta (w_{\sigma ,m})=v_{\sigma
,m}^1$.
\end{quote}
\noindent (The notation $\Theta
_n^{\beta _1,\beta _2} $ is defined before (II.2).)

$\diamondsuit (E_0)$ gives us for each $\sigma \in E_0$ a ``prediction'' $%
\Upsilon (\sigma )=\left\langle \Theta _n^\sigma :n\in \omega \right\rangle $
of an $\omega $-sequence of finite subsets of $T_\sigma $ such that $%
\bigcup_{n\in \omega }\Theta _n^\sigma $ is an antichain in $T_\sigma $. The
proof that $G^0$ and $G^T$ are not isomorphic will then work as before.

\medskip\ 

{\bf (VII)} The next step is to describe in detail the inductive
construction of the data satisfying the properties given above. Our
construction is by induction on the elements of $E$. At stage $\delta \in E$
we will define $G_\alpha ^0$ and $G_\alpha ^\delta $ for any $\alpha \leq
\delta +1$ for which they are not already defined. We will have already
defined $G_\nu ^0$ for $\nu \leq \sup \{\delta ^{\prime }+1:\delta ^{\prime
}\in E$, $\delta ^{\prime }<\delta \}$. By following the prescriptions in
(*) and (***), we can assume that $G_\nu ^0$ is defined for all $\nu \leq
\delta $.

Let $\gamma =\sup \{\delta ^{\prime }+1:\delta ^{\prime }\in E\cap \delta $, 
$T_\delta  \rest \delta ^{\prime }=T_{\delta ^{\prime }}\}$. Then we need to
define $G_\alpha ^\delta $ for $\gamma <\alpha \leq \delta +1$. We need to
do this is such a way that we are able to define the partial isomorphisms $%
f_\alpha ^\delta $. We shall leave the details of the latter to the next
section and describe the construction of the groups here. There are two
cases to consider.

\smallskip\ 

{\it Case 1: }$\gamma =\delta \in E${\it .} Then $G_\delta ^\delta $ is
already defined. If $\delta \in E_0$, follow the prescription in (***) and
(****). If $\delta \in E_1$, $\diamondsuit (E_1)$ gives us an isomorphism $%
h:G_\delta ^0\rightarrow G_\delta ^\delta $; the construction of $G_{\delta
+1}^0$ and $G_\delta ^\delta $ is essentially the same as in the previous
Theorem (Case 4 of (III)); in particular, if (II.2) holds, we use an
antichain $A_{\beta ,\delta }^\delta =\{t:t$ is $<_\delta $-minimal in $%
\delta \setminus \beta \}$; $G_{\delta +1}^0$ is generated by $G_\delta
^0\cup \{z_{\delta ,n}^0:n\in \omega \}$ subject to relations $p_nz_{\delta
,n+1}^0=z_{\delta ,n}^0+k_{\delta ,n}^0$ (which keep $h$ from extending) and 
$G_{\delta +1}^\delta $ is generated by $G_\delta ^\delta \cup \{z_{\delta
,n}^\delta : n\in \omega \}$ subject to relations $p_nz_{\delta ,n+1}^\delta
=z_{\delta ,n}^\delta +k_{\delta ,n}^\delta $ (where $k_{\delta ,n}^\delta $
= $f_{\Theta _n}^\delta (k_{\delta ,n}^0)$).

For the purposes of later stages of the construction we also define, for any 
$\delta _1>\delta $ such that $\delta _1\in E$ and $T_{\delta _1} \rest \delta
\neq T_\delta $, elements $k_{\delta ,n}^{\delta _1}\in G_\delta ^{\delta
_1} $. We know that $k_{\delta ,n}^0$ has the form $x_{\mu _n,j_n}^0\pm \xi
_{\delta ,n}$ where $\xi _{\delta ,n}$ is either $0$, $x_{\sigma ,j}^0$, or $%
w_{\sigma ,j}$ for some $\sigma $, $j$ (cf. (IV.1)). In case $\xi _{\delta
,n}$ is $0$, let $k_{\delta ,n}^{\delta _1}=x_{\mu _n,j_n}^1$; in case $\xi
_{\delta ,n}=x_{\sigma ,j}^0$ , let $k_{\delta ,n}^{\delta _1}=x_{\mu
_n,j_n}^1\pm x_{\sigma ,j}^1$. Finally, if $\xi _{\delta ,n}=w_{\sigma ,j}$,
let $k_{\delta ,n}^{\delta _1}=x_{\mu _n,j_n}^1\pm \xi _{\delta ,n}^{\prime
} $ where%
$$
\xi _{\delta ,n}^{\prime }=\left\{ 
\begin{array}{ll}
w_{\sigma ,j}^1 & \hbox{ if }T_{\delta _1} \rest \sigma \neq T_\sigma \\ v_{\sigma
,j}^1 & \hbox{if }T_{\delta _1} \rest \sigma =T_\sigma 
\end{array}
\right. 
$$
and $w_{\sigma ,j}^1=2u_{\sigma ,j+1}^1-u_{\sigma ,j}^1$. We will be able to
show (in the next section) the following:

\begin{quote}
(VII.1) for any branch $B$ in $T_{\delta _1} \rest \delta $ with $\delta =\sup
\{t+1:t\in B\}$, $g_B=\cup \{f_\alpha ^{\delta _1}:\alpha \in B\}$ is such
that for almost all $n$, $g_B(k_{\delta ,n}^0)=k_{\delta ,n}^{\delta _1}$.
\end{quote}

\noindent (This is evidence of what, in view of (IV.3), will enable us to
extend functions.)

{\it Case 2: $\gamma <\delta $. }We need to define $G_\alpha ^\delta $ for $%
\gamma +1\leq \alpha \leq \delta +1$ by induction on $\alpha $. If we have
defined $G_\alpha ^\delta $ for $\alpha \leq \rho <\delta $, and $\rho $
does not belong to $E_1$, we follow the prescription in (**) or (****). If $%
\rho \in E_1$, then $T_\delta  \rest \rho \neq T_\rho $ (by definition of $\gamma $%
). By induction $G_{\rho +1}^0$ is constructed as in Case 1 and we have $%
k_{\rho ,n}^\delta $ as there (with $\delta $ playing the role of $\delta _1$
and $\rho $ playing the role of $\delta $). In particular, $G_{\rho +1}^0$
is generated by $G_\rho ^0\cup \{z_{\rho ,n}^0:n\in \omega \}$ subject to
relations $p_nz_{\rho ,n+1}^0=z_{\rho ,n}^0+k_{\rho ,n}^0$. We define $%
G_{\rho +1}^\delta $ to be generated by $G_\rho ^\delta \cup \{z_{\rho
,n}^\delta :n\in \omega \}$ subject to relations $p_nz_{\rho ,n+1}^\delta
=z_{\rho ,n}^\delta +k_{\rho ,n}^\delta $. Finally, we define $G_{\delta
+1}^\delta $ as in Case 1.

\smallskip\ 

The definition of the $W_\alpha ^\delta [\Theta ]$ will be as in (III);
specifically, \ $W_{\alpha +1}^\delta [\Theta ]=W_\alpha ^\delta [\Theta ]$
unless $\alpha =\sigma \in E_0$, $T_\delta  \rest \sigma =T_\sigma $ and $\Theta
=\Theta _n^\sigma $ for some $n$, in which case $W_{\sigma +1}^\delta
[\Theta ]=W_\sigma ^\delta [\Theta ]\cup \{w_{\sigma ,n}\}$.

\medskip\ 

{\bf (VIII)}\ We have defined the groups and the sets $W_\alpha ^\delta
[\Theta ]$; the last step is to show that the partial isomorphisms $f_\nu
^\delta $ can be defined satisfying the conditions in 4.

First let us verify (VII.1). Let $\delta $ and $\delta _1$ be
as in Case 1 of (VII) and suppose $B$ is a branch in $T_{\delta _1} \rest \delta $
with $\delta =\sup \{t+1:t\in B\}$. Then $g_B$ is an isomorphism $:G_\delta
^0\rightarrow G_\delta ^{\delta _1}$ and we want to show that $g_B(k_{\delta
,n}^0)=k_{\delta ,n}^{\delta _1}$ for almost all $n$. Recall that $k_{\delta
,n}^0$ has the form $x_{\mu _n,j_n}^0\pm \xi _{\delta ,n}$ where $\xi
_{\delta ,n}$ is either $0$, $x_{\sigma ,j}^0$, or $w_{\sigma ,j}$ for some $%
\sigma $, $j$; the only case we need to worry about is when $\xi _{\delta
,n}=w_{\sigma ,j}$. 

Let $\mu =\sup \{\alpha <\delta :T_\delta  \rest \alpha =T_{\delta _1} \rest \alpha \}$;
so $\mu < \delta $ and $G_\alpha ^{\delta _1}=G_\alpha ^\delta $ for $%
\alpha \leq \mu $. Suppose that $G_{\delta +1}^0$ and $G_{\delta +1}^\delta $
are defined using $A_{\beta ,\delta }^\delta =\bigcup_{n\in \omega }\Theta
_n^{\beta ,\delta }$ as in Case 1 of (VII) and Case 4 of (III). We consider
several cases. First, suppose that there exists $t\in A_{\beta ,\delta
}^\delta $ with $t\geq \mu $. Then for almost all $n$, $t\in \Theta
_n^{\beta ,\delta }$ and thus if $w_{\sigma ,j}\in W_\delta ^\delta [\Theta
_n^{\beta ,\delta }]$ then $\sigma >t\geq \mu $; hence $T_{\delta _1} \rest \sigma
\neq T_\delta  \rest \sigma $ and it follows from 4(d) that $g_B(k_{\delta
,n}^0)=k_{\delta ,n}^{\delta _1}$. If this case does not hold then $A_{\beta
,\delta }^\delta \subseteq \mu $ so $A_{\beta ,\delta }^\delta =A_{\beta
,\mu }^\delta $ is an antichain in $T_{\delta _1} \rest \mu =T_\delta  \rest \mu $. If
there exists $t\in B$ with $\beta \leq t<\mu $, then there exists $t\in B$
with $t\in A_{\beta ,\delta }^\delta $ and hence $t\in \Theta _n^{\beta
,\delta }$ for almost all $n$; it follows easily that for almost all $n$ $%
g_B(k_{\delta ,n}^0)=k_{\delta ,n}^{\delta _1}$ (considering separately the
cases when $\sigma \leq \mu $ and $\sigma >\mu $). In the remaining case, if 
$\alpha =\inf \{t\in B:t\geq \beta \}$, then $\alpha \geq \mu $ so we have $%
\sup \{t:t<_{\delta _1}\alpha \}\leq \beta <\mu \leq \alpha $ and we have
the desired conclusion by 4(e) --- again distinguishing between the cases
when $\sigma \leq \mu $ and $\sigma >\mu $. This completes the proof of
(VII.1).

\smallskip\ 

Now we need to verify the analog of (III.1). Letting $\delta $ and $\gamma $
be as in (VII), we need to define $f_\alpha ^\delta $ for $\gamma \leq
\alpha \leq \delta $. Fix $\alpha $ and let $B=\{t<\gamma :t<_\delta \alpha
\}$ and $g_B=\cup \{f_t^\delta :t\in B\}$. We can suppose that $\alpha $ is $%
<_\delta $-minimal among elements of $\{\beta :\gamma \leq \beta \leq \alpha
\}$.

We will first define an extension of $g_B$ to a partial isomorphism $\tilde
g_B$ whose domain is 
$$
\limfunc{dom}(g_B)+\left\langle 
\begin{array}{c}
\{x_{\nu ,j}^0:\nu 
\hbox{ }\leq \alpha \hbox{, }j=0,1\}\cup \\ \{u_{\sigma ,n}^0:\sigma \in
E_0\cap (\alpha +1) 
\hbox{, }n\in \omega \}\cup \\ \{v_{\sigma ,n}^0:\sigma \in E_0\cap (\alpha
+1)\hbox{, }n\in \omega \} 
\end{array}
\right\rangle 
$$

Using an enumeration in an $\omega $-sequence of $Y_0\cup Y_1$ where%
$$
Y_0=\{\sigma \in E_0:\sup B\leq \sigma <\gamma \hbox{ and }T_\delta  \rest \sigma
=T_\sigma \} 
$$
and%
$$
Y_1=\{\left\langle \beta _1,\beta _2\right\rangle :\sup B\leq \beta _1<\beta
_2\leq \alpha \} 
$$
we can define $\tilde g_B$ such that

\begin{quote}
(c$^{\prime }$) for all $\nu \leq \alpha $, $j\in \{0,1\}$ $\tilde g_B(x_{\nu
,j}^0)=x_{\nu ,j}^1$; moreover, if $w_{\sigma ,n}\in W_{\gamma +1}^\delta
[\Theta ]$ and $\Theta \cap B\neq \emptyset $, then $\tilde g_B(w_{\sigma
,n})=v_{\sigma ,n}^1$;

(d$^{\prime }$) if $\sigma \in E_0\cap \alpha + 2 $, then $\tilde g_B(u_{\sigma
,n}^0)=u_{\sigma ,n}^1$ for almost all $n$, and if $T_\delta  \rest \sigma \neq
T_\sigma $, then $\tilde g_B(u_{\sigma ,n}^0)=u_{\sigma ,n}^1$ for all $n$;
and

(e$^{\prime }$) for all pairs $\beta _1$, $\beta _2$ with $\sup B\leq \beta
_1<\beta _2\leq \alpha $, it is the case for almost all $n\in \omega $ that
for all $w_{\sigma ,m}\in W_{\alpha +1}^\delta [\Theta _n^{\beta _1,\beta
_2}]$ we have $\tilde g_B(w_{\sigma ,m})=v_{\sigma ,m}^1$.
\end{quote}

\smallskip\ 

Now $\tilde g_B(k_{\rho ,n}^0)\,$ is defined for all $\rho \in E_1$ with $%
\rho \leq \alpha $. We need to define $f^\delta _ \alpha (z_{\rho ,n}^\delta
) $ for all such $\rho \geq \sup B$. In view of (IV.3), we can do this
provided that $\tilde g_B(k_{\rho ,n}^0)=k_{\rho ,n}^\delta $ for almost all 
$n\in \omega $. We consider separately the cases: $T_\delta  \rest \rho =T_\rho $;
and $T_\delta  \rest \rho \neq T_\rho $. The first case is as in (IV); the last is
as in the proof of (VII.1) (with $\delta $ playing the role of $\delta _1$, $%
\rho $ playing the role of $\delta $ and using (d$^{\prime }$) and (e$%
^{\prime }$)).

This completes the proof of  Theorem 8.

\end{document}